\documentclass[11pt]{amsart}

\usepackage{amsfonts, enumerate} 
\usepackage{amssymb}
\usepackage[final]{graphicx, graphics,stmaryrd} 
\input xy 
\xyoption{all} 
\usepackage[margin=3cm]{geometry} 
\usepackage{amscd}
\usepackage[all]{xy}
\usepackage[breaklinks,bookmarksopen,bookmarksnumbered]{hyperref}

\usepackage[english]{babel}
\usepackage[utf8]{inputenc}
\usepackage{amsthm}
\usepackage{graphics}
\usepackage{amsmath}
\usepackage{amstext}
\usepackage{engrec}
\usepackage{rotating}
\usepackage[safe,extra]{tipa}
\usepackage{multirow}
\usepackage{enumitem}
\usepackage{bm}

\newtheorem{teo}{Theorem}[section]
\newtheorem{coro}[teo]{Corollary}
\newtheorem{prop}[teo]{Proposition}
\newtheorem{lemma}[teo]{Lemma}

\newtheorem*{teon}{Theorem}

\theoremstyle{definition}
\newtheorem{defin}[teo]{Definition}
\newtheorem{rmk}[teo]{Remark}
\newtheorem{fact}[teo]{Fact}
\newtheorem*{ter}{Terminology}

\newcommand{\mb}{\mathbb}
\newcommand{\mc}{\mathcal}
\newcommand{\msf}{\mathsf}

\def\R{\mathcal R}
\def\C{\mathbb C}

\def\Z{\mathbb Z}

\def\Q{\mathbb Q}
\def\p{\mathbb P}

\def\E{\mathbb E}
\def\M{\mathcal{M}}
\def\A{\mathbb A}

\def\mm{\bm{\mathsf{M}}}
\def\rr{\bm{\mathsf{R}}}

\newcommand{\ssC}{\mathsf{C}}
\newcommand{\ssL}{\mathsf{L}}

\newcommand{\m}{\mbox}

\newcommand{\cor}{\textit}

\newcommand{\fine}{\qed\newline}
\newcommand{\xx}{\otimes}

\newcommand{\beginCD}{\begin{equation*}\begin{CD}}
\newcommand{\enCD}{\end{CD}\end{equation*}}
\newcommand{\mdeg}{\underline\deg}

\DeclareMathOperator{\age}{age}

\DeclareMathOperator{\defo}{Def}
\DeclareMathOperator{\Aut}{Aut}
\DeclareMathOperator{\diag}{Diag}

\DeclareMathOperator{\ord}{ord}
\DeclareMathOperator{\im}{Im}
\DeclareMathOperator{\id}{id}
\DeclareMathOperator{\sing}{Sing}

\DeclareMathOperator{\Ker}{Ker}

\DeclareMathOperator{\GL}{GL}

\DeclareMathOperator{\QR}{QR}
\DeclareMathOperator{\Pic}{Pic}

\DeclareMathOperator{\cut}{cut}
\DeclareMathOperator{\Rep}{Rep}

\DeclareMathOperator{\codim}{Codim}

\DeclareMathOperator{\nc}{nc}

\newcommand{\mmu}{{\pmb \mu}}

\DeclareMathOperator{\sep}{sep}

\def\s{\mc S}

\DeclareMathOperator{\cord}{c-ord}

\title{Singularities of moduli of curves with a universal root}
\author{Mattia Galeotti}

\begin{document}

\begin{abstract}
In a series of recent papers, Chiodo, Farkas and Ludwig carry out a
deep analysis of the singular locus of the moduli space of stable (twisted) curves
with an $\ell$-torsion line bundle. They show that for $\ell\leq 6$ and $\ell\neq 5$
pluricanonical forms extend over any desingularization.
This opens the way to a computation of the Kodaira dimension without desingularizing,
as done by Farkas and Ludwig for $\ell=2$,
and by Chiodo, Eisenbud, Farkas and Schreyer
for $\ell=3$.
Here we treat roots of line bundles on the universal curve systematically:
we consider the moduli space of curves
$C$ with a line bundle $L$ such that $L^{\xx\ell}\cong\omega_C^{\xx k}$.
New loci of canonical and non-canonical singularities appear
for any $k\not\in\ell\Z$ and $\ell>2$, we provide a set of combinatorial tools allowing us
to completely describe the singular locus in terms of dual graphs.
We characterize the locus of non-canonical singularities, and for small
values of $\ell$ we give an explicit description.

\end{abstract}

\maketitle

The moduli space  $\M_g$ of smooth curves of genus $g$,
alongside with its compactification $\overline\M_g$ in the sense of Deligne-Mumford,
is one of the most fascinating and largely studied objects in algebraic geometry.
 It is also the subject of wide open questions; the one that motivates our work here is inspired by Harris and Mumford
\cite{harmum82}, they prove that for $g>23$, $\M_g$ is of general type.
The question on how and precisely where the transition to general type varieties happens
remains open and highly intertwined with several other open questions in this field.
However, a new line of attack seems to stand out clearly.
It appears that the study of the finite covers of moduli of curves and roots
of a universal bundle is more manageable. 
A series of papers by Farkas and Farkas-Verra (see \cite{far00} and \cite{farver10}) provide a complete description for spin structures.
Chiodo, Eisenbud, Farkas and Schreyer show in \cite{cefs13} that for curves with a $3$-torsion bundle,
the transition happens before $g=12$.
The previous results use a compactification of the moduli spaces via twisted curves following from \cite{abravis02,acv03,abrajar03,chio08},
and are based on a preliminary analysis by Chiodo, Farkas and Ludwig of the singular locus 
of the moduli space of curves carrying a torsion line bundle, or level curves:
in particular in \cite{chiofar12} they characterize the singular locus and show that
for $\ell$-torsion bundles with $\ell\leq 6$ and $\ell\neq 5$, pluricanonical
forms could be extended over any desingularization using
a method developed by Harris and Mumford.

The purpose of this paper is to identify and analyze the locus of
singularities for every line bundle
on the universal curve. If $\mc{U}_g$ is the universal curve on $\M_g$,
we can focus on roots of any order of the powers of the relative dualizing sheaf $\omega=\omega_{\mc{U}_g\slash\M_g}$.
Indeed, these are all possible bundles on $\mc{U}_g$, up to pullbacks
from the moduli space $\M_g$, as discussed for example in \cite{mestrarama85, arbacor87}.
We consider the moduli space $\R_{g,\ell}^k$
parametrizing $\ell$th roots of $\omega^{\xx k}$, \cor{i.e.}~triples
$(C,L,\phi)$ where $C$ is a smooth curve, $L$ a line bundle over $C$
and $\phi$ an isomorphism $L^{\xx \ell}\to \omega_C^{\xx k}$.
In the case $k=0$, our analysis matches with Chiodo and Farkas work on
level curves \cite{chiofar12}. Furthermore, we give a more general
description of the singularities for all $k$, in particular on those cases
where it is not possible to extend everywhere the pluricanonical
forms using the known techniques.

To compactify 
$\R_{g,\ell}^k$ we use the theory of stack-theoretic stable curves:
we consider the moduli stack $\overline\rr_{g,\ell}^k$ sending
a scheme $S$ to a triple $(\ssC\to S,\ssL,\phi)$ where
$\ssC$ is a stack-theoretic curve with a stable curve $C$ as coarse space and
possibly non-trivial stabilizers on its nodes, $\ssL$ is a line bundle over $\ssC$
and $\phi$ an isomorphism $\ssL^{\xx \ell}\to \omega_\ssC^{\xx k}$.
This moduli is represented by a Deligne-Mumford stack.
The open dense substack $\rr_{g,\ell}^k$, which admits only smooth curves~$\ssC$, has
$\R_{g,\ell}^k$ as coarse space.
The stack structure extends the covering $\R_{g,\ell}^k\to \M_g$
on all $\overline\M_g$.\newline

In this paper we focus on the singular locus of $\overline\R_{g,\ell}^k$,
called $\sing\overline\R_{g,\ell}^k$, and
on the locus of non-canonical singularities, also called non-canonical locus and noted $\sing^{\nc}\overline\R_{g,\ell}^k$. 
The singular locus of $\overline\M_g$
is characterized in \cite{harmum82} by introducing the concept of 
elliptic tail quasireflection
(see Definition \ref{def_eti}): $[C]\in\overline\M_g$ is in $\sing\overline\M_g$ if and only if there exists
$\alpha\in\Aut(C)$ not generated by elliptic tail quasireflections. This locus 
admits a natural lift to~$\overline\R_{g,\ell}^k$ by introducing the group $\Aut'(C)$
of automorphisms on $C$ which preserve the root structure. However, the stack
structure of a curve comes equipped with new \cor{ghost} automorphisms, \cor{i.e.}~automorphisms
acting trivially on the coarse space of the curve, but possibly non-trivially
on the additional structure. Because of such ghosts, there is a new
locus of singularities which can be studied
using new tools in graph theory.

First, to each point $[\ssC,\ssL,\phi]\in\overline\R_{g,\ell}^k$
we attach its dual graph $\Gamma(\ssC)$: its vertex set is,
as usual, the set of irreducible components of $\ssC$, and its edge set
is the set of nodes of $C$. Moreover, we introduce a multiplicity index
$M$ (see Definition \ref{mul_index}), which is a function on the edges of the dual graph and depends
on the line bundle $\ssL$ (it maps each node with a privileged branch to the character $\ssL|_{\msf n}$).
Finally, for any prime $p$ in the factorization of
$\ell$, if $e_p$ is the $p$-adic valuation power of $\ell$,
we note $\Gamma_p(\ssC)$ the graph obtained
contracting those edges of the dual graph where $M$ is divisible by $p^{e_p}$.
A graph is tree-like if it has no circuits except for loops.
Generalizing previous results from \cite{lud10} and \cite{chiofar12},
we prove the following (see Theorem \ref{smooth1}).
\begin{teon}
For any $g\geq 4$ and $\ell$ positive integer, the point $[\ssC,\ssL,\phi]\in\overline\R_{g,\ell}^k$ is smooth
if and only if the two following conditions are satisfied:
\begin{enumerate}[label=\roman{*}., ref=(\roman{*})]
\item  the image $\Aut'(C)$ of the coarsening morphism $\underline\Aut(\ssC,\ssL,\phi)\to\Aut(C)$ is generated by elliptic tail quasireflections;
\item the dual graph contractions $\Gamma(\ssC)\to\Gamma_p(\ssC)$ yield a tree-like
graph for
all prime numbers $p$ dividing $\ell$.
\end{enumerate}
\end{teon}

This is a generalization of Chiodo and Farkas \cite[Theorem 2.28]{chiofar12}.
In the case $k=0$ treated by them, the dual graph contractions must be 
\cor{bouquet} graphs, \cor{i.e.}~graphs with only one vertex.
We spell out this generalization on Remark \ref{rmk_teo}.

The role of dual graph structure in characterizing the geometry of $\overline\R_{g,\ell}^k$
goes even deeper.
Consider the contraction $\Gamma(\ssC)\to\Gamma_0(\ssC)$ of the edges
where $M$ vanishes.
The graph $\Gamma_0(\ssC)$ and the multiplicity index $M$
completely describe the group of ghost automorphisms, and as a consequence
they carry all the information about the new singularities. 
We have that the non-canonical locus is of the form
$$\sing^{\nc}\overline\R_{g,\ell}^k=T_{g,\ell}^k\cup J_{g,\ell}^k,$$
where $T_{g,\ell}^k\subset\pi^{-1}\sing^{\nc}\overline\M_g$ is the analogue
of old non-canonical singularities, and $J_{g,\ell}^k$ is the 
new part of the non-canonical locus. 
For the $T$-locus
the generalization of the technique developed in \cite{harmum82}
allows the extension of pluricanonical forms. Unfortunately, it is
unknown if this extension is possible over the $J$-locus.
This is the main obstacle toward the analysis of birational geometry $\overline\R_{g,\ell}^k$.
One of the aims of this paper is precisely to describe $J_{g,\ell}^k$ further.
In Section \ref{subs_strata} we introduce a
stratification of the boundary $\overline\R_{g,\ell}^k\backslash\R_{g,\ell}^k$
labelled by decorated graphs, \cor{i.e.}~pairs $(\Gamma_0(\ssC),M)$:
we will decompose the $J$-locus using this stratification.

In the last section we work out the cases $\ell=2,\ 3,\ 5$ and $7$.
Using the age criterion, we develop a machinery to study the ghosts
group through graph properties, we conclude in particular
that 
in many cases the $J$-locus is the union of strata labelled only by \cor{vine} graphs,
\cor{i.e.}~graphs with only two vertices.
The first exception comes out with $\ell=7$, when a $3$-vertices graph stratum is
needed to complete the description.

\section*{Acknowledgements}
This paper is part of my Ph.D. research, 
which I am doing at the
Institut de Mathématiques de Jussieu-Paris Rive Gauce,
under the supervision of Alessandro Chiodo.
I am particularly grateful to him for his support, for the motivating discussions and for his patience.
I am also very grateful to the anonymous referee for his careful reading and his comments.

\section{Preliminaries on curves and their dual graphs}

In this paper $\ell$ is always a positive integer. Every curve is an algebraic curve of genus $g\geq 0$ over an algebraically
closed field~$k$ whose characteristic is $0$ or prime to $\ell$. By $S$-curve we mean a family of curves
on a scheme $S$, \cor{i.e.}~a flat morphisms $C\to S$ such that every geometric fibre is a curve of given genus.

Given a line bundle $F$ on a curve $C$, an $\ell$\cor{th root} of $F$ on $C$ is a pair 
$(L;\phi)$,
where $L$ is a line bundle on $C$ and
$\phi$ is an isomorphism of line bundles
$ L^{\otimes\ell}\to F$
is an isomorphism. 
A morphism between two $\ell$th roots $(L,\phi),\ (L',\phi')$, 
is a morphism $\rho\colon L\to L'$ of line bundles on $C$ such that
$\phi'=\phi\circ\rho^{\xx\ell}$.

If $C$ is a smooth curve and $F$ is a line bundle on $C$
and $\ell$ divides $\deg F$, $F$ has at least one $\ell$th root, and in this case the number of $\ell$th
roots is exactly $\ell^{2g}$. This property fails for
general singular curves. 
Replacing scheme theoretical curves by twisted curves, that we define
just below, we still get the ``right'' number of $\ell$th roots.

In the last part of this section  we will introduce
some basic tools on graph theory.

\subsection{Twisted curves and $\ell$th roots}

We recall the definition of twisted curve.

\begin{defin}[twisted curve]
A twisted curve $\ssC$ is a Deligne-Mumford stack whose coarse space is a curve with nodal singularities,
such that its smooth locus is represented
by a scheme, and the local pictures of the singularities are given by
$[\{xy=0\}\slash{\pmb \mu}_r]$, $r$ positive integer, with any primitive $r$th root of unity $\zeta$ acting as $\zeta\cdot (x,y)=(\zeta x,\zeta^{-1}y)$.
In this case we say that the node has a $\mmu_r$ stabilizer.
\end{defin}

If $\msf{n}$ is a node of $\ssC$ with non-trivial stabilizer
$\mmu_r$, the local picture of the curve at $\msf{n}$ is $xy=0$ with the action of $\mmu_r$
described above. Given a line bundle $\msf{F}\to\ssC$ at $\msf{n}$,
its local picture at the node is $\mathbb{A}^1\times \{xy=0\}\to\{xy=0\}$
with any primitive root $\zeta\in\mmu_r$ acting as
\begin{equation}\label{node_twist}
\zeta\cdot (t,x,y)=(\zeta^mt,\ \zeta x,\ \zeta^{-1} y),\ \ \m{with}\ \ m\in\Z\slash r,
\end{equation}
on $\mathbb{A}^1\times \{xy=0\}$. The index $m\in \Z\slash r$, called
\cor{local multiplicity}, is uniquely determined when we assign a privileged choice
of a branch of $\mathsf{n}$ where $\zeta$ acts as $x\mapsto \zeta x $. If we switch
the privileged branch, the local multiplicity changes to its opposite $-m$.

A line bundle $\ssL$ on $\ssC$ is \emph{faithful} if the associated morphism
$\ssC\to\msf{B}\mathbb{G}_m$ is representable. As pointed out in \cite[Remark 1.4]{chiofar12},
a line bundle on a twisted curve is faithful if and only if the local multiplicity index $m$ at every
node $\msf{n}$ is coprime with the order $r$ of the local stabilizer.

\begin{defin}
An $\ell$th root of the line bundle $\msf{F}$ on the twisted curve $\ssC$, is a pair
$(\ssL,\phi)$,
where $\ssL$ is a faithful line bundle on $\ssC$
and $\phi\colon \ssL^{\otimes \ell}\to \mathsf{F}$ 
is an isomorphism between $\ssL^{\otimes\ell}$ and $\msf{F}$.
\end{defin}

Fixing a primitive $r$th
root of the unity $\zeta$ will not affect our results and will permit
us to identify (non-canonically) $\mmu_r$ and $\Z\slash r$.
In what follows we chose the primitive root $\xi_r=\exp(2i\pi\slash r)$ for every node.

\subsection{The dual graph}

If a curve is not smooth, it may have several irreducible components.
The information about the relative crossing of this components is
encoded in the dual graph.

\begin{defin}[dual graph]
Given a curve $C$ (scheme theoretic or twisted), the dual graph $\Gamma(C)$ has
the set of irreducible components of $C$ as vertex set, 
and the set of nodes of $C$ as edge set $E$. 
The edge associated to the node $\msf{n}$
links the vertices associated to the components where the branches of $\msf{n}$ lie.
\end{defin}

We introduce some graph theory which will be
important to study the moduli spaces of twisted curves
equipped with a root.

\subsubsection{Cochains over a graph}\label{coch_sect}
Consider a connected graph $\Gamma$ with vertex set $V$ and edge set~$E$, we call
\cor{loop} an edge that starts and ends on the same vertex, we call
\cor{separating} an edge $e$ such that the graph with vertex set $V$ and
edge set $E\backslash \{ e\}$ is disconnected. We note by $E_{\sep}$ the set
of separating edges. We note by $\E$ the set of oriented
edges: the elements of this set are edges in $E$ equipped with an orientation.
In particular for every edge $e\in\E$ we note
$e_+$ the head vertex and $e_-$ the tail, and there is a $2$-to-$1$ projection $\E\to E$. We also introduce a conjugation 
in $\E$, such that for each $e\in \E$,  the conjugated edge $\bar e$ is obtained
by reversing the orientation, in particular
$(\bar e)_+=e_-$.

\begin{rmk}\label{rmk_or}
In the dual graph setting, where the vertices are the components of a curve $C$
and the edges are the nodes, the oriented edge set $\E$ could be seen as the set of
branches at the curve nodes. Indeed, every edge $e$, equipped with an orientation,
is bijectively associated to the branch it is pointing at.
\end{rmk}

\begin{defin}
The group of $0$-cochains is the group of $\Z\slash\ell$-valued functions on $V$
$$C^0(\Gamma;\Z\slash\ell):=\left\{a\colon V\to \Z\right\}=\bigoplus_{v\in V}\Z\slash\ell.$$
The group of $1$-cochains is the group of antisymmetric functions on $\E$
$$C^1(\Gamma;\Z\slash\ell):=\left\{b\colon \E\to\Z\slash\ell|\ b(\bar e)=-b(e)\right\}.$$
After assigning an orientation to every edge $e\in E$, we 
may identify $C^1(\Gamma)=\bigoplus_{e\in E}\Z\slash\ell$, but
we prefer working without any prescribed choice of orientation.
\end{defin}

These spaces are equipped with two non-degenerate bilinear $\Z\slash\ell$-valued forms
$$\langle a_1,a_2\rangle:=\sum_{v\in V}a_1(v)a_2(v)\ \ \ \ \langle b_1,b_2\rangle:=\frac{1}{2}\sum_{e\in \E}b_1(e)b_2(e)$$
for all $a_1,a_2\in C^0(\Gamma)$ and $b_1,b_2\in C^1(\Gamma)$.
Consider
the operator $\partial\colon C^1(\Gamma)\to C^0(\Gamma)$ such that
$$\partial b(v):=\sum_{e_+=v}b(e),\ \ \forall b\in C^1(\Gamma)\ \forall v\in V.$$
\begin{defin}
We introduce $G(\Gamma;\Z\slash\ell):=\left(\ker\partial\right)^\perp$,
the orthogonal complement with respect to the form just introduced.
\end{defin}
We also introduce the 
exterior differential $\delta\colon C^0(\Gamma)\to C^1(\Gamma)$ such that
$$\delta a(e):=a(e_+)-a(e_-),\ \ \forall a\in C^0(\Gamma)\ \forall e\in \E.$$
\begin{prop}[see {\cite[p.9]{chiofar12}}]\label{prop_G}
The adjoint of $\partial$ is $\delta$, and
$G(\Gamma;\Z\slash\ell)=\im\delta$.\fine
\end{prop}

The exterior differential fits into an useful exact sequence.
$$
0\to\Z\slash\ell\xrightarrow{i} C^0(\Gamma;\Z\slash\ell)\xrightarrow{\delta} C^1(\Gamma;\Z\slash\ell).
$$
Where the injection $i$ sends $m\in\Z$ on the cochain constantly equal to $m$.
This sequence gives the dimension of $\im\delta$.
\begin{equation}
\im\delta\cong (\Z\slash\ell)^{\#V-1}.
\end{equation}

\subsubsection{Construction of a basis for $C^1(\Gamma)$} 

\begin{defin}[cuts, paths and circuits]
A cut is a $1$-cochain $b\colon\E\to \Z$ such that there exists a non-empty subset $W\subset V$
and the values of $b$ are determined in the following way:
 $b(e)=1$ if the head of $e$ is in $W$ and the tail in $V\backslash W$, $b(e)=-1$
if the head is in $V\backslash W$ and the tail in $W$, $b(e)=0$ elsewhere.
The support of $b$ in $\E$ is sometimes called \cor{cut-set} of $b$.
 
A path in a graph $\Gamma$ is a sequence $e_1, e_2,\dots, e_k$ of edges in $\E$
overlying $k$ distinct non-oriented edges in $E$, and
such that the head of $e_i$ is the tail of $e_{i+1}$ for all $i=1,\dots, k$.

A circuit is a closed path, \cor{i.e.}~a path $e_1,\dots, e_k$ such that
the head of $e_k$ is the tail of $e_1$.
We often refer to a circuit by referring to its characteristic function, \cor{i.e.}~the
 cochain $b$ such that $b(e_i)=1$ for all $i$, $b(\bar e_i)=-1$ and $b(e)=0$ if $e$
is not on the circuit.
\end{defin}
\begin{rmk}
Any cut $b$ is an element of $\im\delta$. Indeed, $b=\delta a$ if $a$ is
the characteristic function of $W$.
\end{rmk}
Moreover,
the group $\im\delta$ is sometimes called \cor{cuts space} because
it is generated by cuts (see Proposition \ref{cut_gen}). 
We also have that
 $\Ker\partial$ is generated by circuits and this implies the following fact.

\begin{prop}\label{prop_imdelta}
A $1$-cochain $b$ is in $G(\Gamma;\Z\slash\ell)=\im\delta$ if and only if, for every circuit $K=(e_1,\dots, e_k)$ in $\E$, we have
$$b(K) =\sum_{i=1}^k b(e_i)=0.$$
\end{prop}

\begin{ter}[graphs tree-like and trees]
A tree-like graph is a connected graph whose only circuits are loops.
A tree is a graph that does not contain any circuit.
\end{ter}
\begin{rmk}\label{ineq_graph}
For every connected graph $\Gamma$, 
the first Betti number $b_1(\Gamma)=\#E-\#V+1$ is the dimension
rank of the homology group $H_1(\Gamma;\Z)$. 
Note that, $b_1$ being positive,
$$\#E\geq \#V-1$$
with equality if and only if $\Gamma$ is a tree.
\end{rmk}
For every connected graph $\Gamma$ with vertex set $V$ and
edge set $E$, we can choose a connected subgraph $T$ with the same vertex set
and edge set $E_T\subset E$ such that $T$ is a tree. The graph
$T$ is a \cor{spanning tree} of $\Gamma$ and has several
interesting properties. We call $\mb{E}_T$ the set of oriented edges of the spanning
tree $T$. Here we notice that E contains a distinguished subset of edges $E_{\sep}$ whose size is \emph{smaller} than V-1
\begin{lemma}\label{ineq_graph2}
If $E_{\sep}\subset E$ is the set of edges in $\Gamma$ that are
separating, then $\#E_{\sep}\leq \#V-1$ with equality if and only if $\Gamma$
is  tree-like.
\end{lemma}
\proof
If $T$ is a spanning tree for $\Gamma$ and $E_T$ its edge set, then
$E_{\sep}\subset E_T$. Indeed, an edge $e\in E_{\sep}$ is the only path between its two extremities, therefore,
since $T$ is connected, $e$ must be in $E_T$.
Thus $\#E_{\sep}\leq\#E_T=\#V-1$, with equality if and only if all the edges
of $\Gamma$ are loops or separating edges, \cor{i.e.}~if $\Gamma$ is a tree-like graph.\fine

\begin{lemma}[{see \cite[Lemma 5.1]{biggs93}}]
If $T$ is a spanning tree for $\Gamma$,
for every oriented edge $e\in\mb{E}_T$ there is a unique
cut $b\in C^1(\Gamma;\Z\slash\ell)$ such that $b(e)=1$ and
the elements of the support of $b$ other than $e$ and $\bar e$ are all
in $\E\backslash \mb{E}_T$.\fine
\end{lemma}
We call $\cut_\Gamma(e;T)$
the unique cut $b$ resulting from the lemma for all $e\in\E_T$.

\begin{prop}\label{cut_gen}
If $T$ is a spanning tree for $\Gamma$, then the elements $\cut_\Gamma(e;T)$,
with $e$ varying on $E_T$,
form a basis of $\im(\delta)$.
\end{prop}
For a proof of this and a deeper
analysis of the image of $\delta$, see \cite[Chap.5]{biggs93}.
\newline

\subsubsection{Definition and basic properties of edge contraction}\label{edgecontraction}
Another tool in graph theory is edge contraction.
If we have a graph $\Gamma$ with vertex set $V$ and
edge set $E$, we can choose a subset $F\subset E$. Contracting
edges in $F$ means taking the graph $\Gamma_0$ such that:
\begin{enumerate}
\item the edge set is $E_0:=E\backslash F$;
\item given the relation in $V$, $v\sim w$ if $v$ and $w$ are linked by an edge $e\in F$, the
vertex set is $V_0:=V\slash\sim$.
\end{enumerate}
The construction of $\Gamma_0$ follows and we have a natural
morphism $\Gamma\to\Gamma_0$ called \cor{contraction} of $F$.
Edge contraction will be useful, in particular
we will consider the image of the exterior differential $\delta$
and its restriction over contraction of a given graph.
If $\Gamma_0$ is a contraction of $\Gamma$, 
then $E(\Gamma_0)$ is canonically a subset of $E(\Gamma)$.
As a consequence, cochains over $\Gamma_0$
are cochains over $\Gamma$ with the additional condition that
the values on $E(\Gamma)\backslash E(\Gamma_0)$ are all $0$. Then
we have a natural immersion
$$C^i(\Gamma_0;\Z\slash\ell)\hookrightarrow C^i(\Gamma;\Z\slash\ell).$$
 Consider the two operators
$$\delta\colon C^0(\Gamma;\Z\slash\ell)\to C^1(\Gamma;\Z\slash\ell)\ \ \ \m{and}\ \ \ \delta_0\colon C^0(\Gamma_0;\Z\slash\ell)\to C^1(\Gamma_0;\Z\slash\ell).$$
Clearly $\delta_0$ is the restriction of $\delta$ on $C^0(\Gamma_0;\Z)$.
From this observation we have the following.
\begin{prop}\label{prop_imcap}
$$\im\delta_0=C^1(\Gamma_0;\Z)\cap \im\delta.$$
\end{prop}
\begin{rmk}\label{rmk_sep}
We make an observation about separating edges of a graph
after contraction. Given a graph contraction $\Gamma\to \Gamma_0$,
the separating edges who are not contracted remain separating. Moreover,
an edge that is not separating cannot become separating. 
\end{rmk}

\section{The moduli space of $\ell$th roots of $\omega^{\xx k}$}

We start by considering the moduli stack $\overline\mm_g$.
Its objects are flat morphisms 
 $C\to S$ such that every geometric fiber
 is a stable genus $g$ curve, \cor{i.e.}~a curve
with at most nodal singularities and an ample canonical line bundle.
The $\overline\mm_g$ coarse space is the moduli space of stable curves $\overline\M_g$.
Harris and Mumford proved in \cite{harmum82} that $\overline\M_g$ 
is of general type for $g> 23$, and classified its
singularities.  We observe that the moduli stack $\mm_g$ of smooth curves is an open
dense substack of $\overline\mm_g$. 

Despite the construction of the moduli space of $\ell$th roots over
smooth curves is quite natural, its extension over $\overline\M_g$ requires some
twisted curve machinery. More precisely, we consider the stack $\rr_{g,\ell}^k$
whose objects are smooth curves $C$ with a line bundle $L$ and an isomorphism
$L^{\xx\ell}\cong \omega_C^{\xx k}$. This stack comes with a natural proper forgetful functor
$\rr_{g,\ell}^k\to\mm_g$. To extend this forgetful morphism over $\overline\mm_g$, we will
introduce $\ell$th roots for any stable curve.

\begin{defin}
Given a stable curve $C$ of genus $g$ and a line bundle $F$ on it, an $\ell$th root of $F$ on $C$
is a triple $(\ssC,\ssL,\phi)$ such that $\ssC$ is a twisted curve, $u\colon \ssC\to C$ its coarsening,
$\ssL$ is a faithful line bundle on $\ssC$ and $\phi\colon\ssL^{\xx\ell}\to u^*F$ an isomorphism.
\end{defin}

\begin{rmk}\label{prop_count}
There exists such a root if and only if $\ell$ divides $\deg F$, and in this case the number of $\ell$th roots
is $\ell^{2g}$.
\end{rmk}

\begin{defin}
The stack $\overline\rr_{g,\ell}^k$ has for objects
the triples $(\ssC\to S,\ssL,\phi)$
where $\ssC\to S$ is a twisted $S$-curve whose coarse space is a stable curve of genus $g$,
$\ssL$ is a faithful line bundle over $\ssC$ and
$$\phi\colon \ssL^{\xx \ell}\to \omega_\ssC^{\xx k}$$
is an isomorphism between $\ssL^{\xx \ell}$ and the $k$th power of the canonical bundle $\omega_\ssC$.
\end{defin}
\begin{rmk}
The canonical bundle $\omega_\ssC$ is the pullback of $\omega_C$ via the coarsening morphism.
In other words, the objects of $\overline\rr_{g,\ell}^k$ are stable curves with an $\ell$th root of their
canonical bundle.
\end{rmk}
\begin{rmk}
Because of Remark \ref{prop_count}, we need that $\ell$ divides $\deg\omega_C^{\xx k}=k\cdot(2g-2)$
to make $\rr_{g,\ell}^k$ non-empty.
\end{rmk}

A triple $(\ssC,\ssL,\phi)$ as above will be improperly called a \cor{rooted curve} when there is no risk of confusion.
The stack $\rr_{g,\ell}^k$ is an open dense substack of $\overline\rr_{g,\ell}^k$.
Moreover, we note $\pi\colon\overline\rr_{g,\ell}^k\to\overline \mm_g$ the extended forgetful functor sending
$(\ssC\to S,\ssL,\phi)$ to $C\to S$ where $C$ is the coarse space of~$\ssC$.
By Remark \ref{prop_count}, every $\pi$ fiber has the same length $\ell^{2g}$.

\subsection{Local structure of the roots moduli}
We consider the stack $\overline\rr_{g,\ell}^k$, with $g\geq 2$ and $0\leq k<\ell$, of
rooted curves of genus $g$. We denote by $\overline\R_{g,\ell}^k$ its associated coarse space.

We are interested in the singularities of this moduli space, in particular
we want a characterization of the singular locus of $\overline\R_{g,\ell}^k$.
If we note $\widetilde \mm_g$ the moduli stack of twisted curves, the forgetful morphism $\overline\rr_{g,\ell}^k\to\widetilde\mm_g$
is étale, as showed in \cite{chio08}. At the same time the composition morphism $\overline\rr_{g,\ell}^k\to\widetilde\mm_g\to\overline\mm_g$
is ramified, as showed by Chiodo and Farkas in \cite{chiofar12}. However, the 
local structure of $\overline\rr_{g,\ell}^k$ is very similar to the
one of $\overline\mm_g$, and we will follow the same approach developed by Chiodo and Farkas
for the case $k=0$.
We denote by $[\ssC,\ssL,\phi]$ the point of $\overline\rr_{g,\ell}^k$ associated to
the rooted curve $(\ssC,\ssL,\phi)$,
then the local picture of the stack $\overline\rr_{g,\ell}^k$ at $[\ssC,\ssL,\phi]$ is
$$\left[\defo(\ssC,\ssL,\phi)\slash\Aut(\ssC,\ssL,\phi)\right],$$
where the universal deformation $\defo(\ssC,\ssL,\phi)$ is a smooth scheme
of dimension $3g-3$, and
$$\Aut(\ssC,\ssL,\phi)=\left\{(\msf{s},\rho)|\ \msf{s}\in\Aut(\ssC)\m{ and }\rho\colon\msf{s}^*\ssL\xrightarrow{\cong}\ssL
\m{ such that }\phi\circ\rho^{\xx\ell}=\msf{s}^*\phi\right\}$$
is the automorphism group of $(\ssC,\ssL,\phi)$. This implies that
the local picture of the moduli space $\overline\R_{g,\ell}^k$
is the classical quotient $\defo(\ssC,\ssL,\phi)\slash\Aut(\ssC,\ssL,\phi)$.
The deformation space
$\defo(\ssC,\ssL,\phi)$ is canonically isomorphic to $\defo(\ssC)$ via the
\'etale forgetful functor $(\ssC,\ssL,\phi)\mapsto \ssC$. Also we
see that the action of $\Aut(\ssC,\ssL,\phi)$ on $\defo(\ssC)$ is not faithful. In particular
the quasi-trivial automorphisms $(\id_\ssC,\zeta)$ with $\zeta\in\mmu_\ell$, whose action scale the fibers,
have trivial action. Thus it becomes natural to consider the group
$$\underline\Aut(\ssC,\ssL,\phi):=\Aut(\ssC,\ssL,\phi)\slash\{(\id_\ssC;\zeta)|\ \zeta\in\mmu_\ell\}=
\left\{\msf{s}\in\Aut(\ssC)|\ \msf{s}^*\ssL\cong\ssL\right\}.$$
The local picture of $\overline\R_{g,\ell}^k$ at $(\ssC,\ssL,\phi)$ could be rewritten as
$$\defo(\ssC)\slash\underline\Aut(\ssC,\ssL,\phi).$$
\begin{defin}
A quasireflection is an element $q\in\GL(\C^m)$ such that its
fixed space is a hyperplane of $\C^m$.
\end{defin}

The smoothness of quotient singularities depends on the 
quotient group to be generated by quasireflections (see \cite{prill67}).
The following fact will permit us to classify moduli space singularities.
\begin{fact}
The scheme theoretic quotient $\defo(\ssC)\slash\underline\Aut(\ssC,\ssL,\phi)$
is smooth if and only if $\underline\Aut(\ssC,\ssL,\phi)$ is spanned by elements
acting as the identity or as quasireflectons.
\end{fact}

Before studying the action of $\underline\Aut(\ssC,\ssL,\phi)$, we point out
the structure of the universal deformation $\defo(\ssC,\ssL,\phi)=\defo(\ssC)$.
\begin{rmk}\label{defo}
We will follow the usual construction developed in \cite{harmum82}. We note $C$
the coarse space of $\ssC$ and
by $N:=\{{\msf n}_1,\dots \msf n_k\}$ the node set of $\ssC$.
If $\{\ssC_1,\dots, \ssC_m\}$ is the set of irreducible components of $\ssC$, we
denote by $\overline \ssC_i$ the normalizations of $\ssC_i$,
by $g_i$ the genus of $\overline C_i$ and
by $\msf{nor}\colon\bigsqcup_{i=1}^m \overline \ssC_i\to \ssC$ the normalization morphism of $\ssC$.
We observe that the curves $\overline\ssC_i$ are scheme theoretic smooth curves.
The divisor $D_i$ on $\overline \ssC_i$ is the preimage of $N$ by the restriction of $\msf{nor}$
on this component. Following \cite{chiofar12} we call $\defo(\ssC,\sing\ssC)$ the universal
deformation of $\ssC$ alongside with its nodes.
Then we have the canonical splitting
$$\defo(\ssC,\sing\ssC)=\bigoplus_{i=1}^m H^1(\overline \ssC_i,T_{\overline \ssC_i}(-D_i)),$$
where $T_{\overline \ssC_i}$ is the tangent bundle to $\overline \ssC_i$
and $H^1(\overline \ssC_i,T_{\overline\ssC_i}(-D_i))=\defo(\overline \ssC_i,D_i)$ parametrizes deformations of
curve $\overline \ssC_i$ with a marking on the points of divisor $D_i$.

Once we mod out $\defo(\ssC,\sing\ssC)$, we have another canonical splitting:
$$\defo(\ssC)\slash\defo(\ssC,\sing\ssC)=\bigoplus_{i=1}^k R_i,$$
where every $R_i$ is unidimensional. The non-canonical choice of a coordinate $t_i$
over $R_i$ corresponds to fixing a smoothing of node $\msf n_i$.
\end{rmk}

Given the coarse space $C$ of $\ssC$, we consider its universal deformation $\defo(C)$
and its universal deformation alongside with its nodes $\defo(C,\sing C)$.
We have another canonical splitting $\defo(C)\slash\defo(C,\sing C)=\bigoplus N_i$,
where the $N_i$ are unidimensional and there are natural coverings $R_i\to N_i$
of order $r(\msf n_i)$ ramified at the origin, where $r(\msf n_i)$
is the order of the stabilizer at $\msf n_i$.

\begin{rmk}\label{gen_nodal}
We point out a classical formula relating the genus of $\ssC$ and the genera of the curves $\ssC_i$.
$$ g(\ssC)=\sum_{i=1}^\nu g_i+b_1(\Gamma)=\sum_{i=1}^\nu g_i+\#E(\Gamma)-\#V(\Gamma)+1,$$
where $\Gamma=\Gamma(\ssC)$ is the dual graph of $\ssC$.
\end{rmk}

The coarsening $\ssC\to C$ induces a group homomorphism
$$\underline\Aut(\ssC,\ssL,\phi)\to \Aut(C).$$
We note the kernel and the image of this morphism by
 $\underline\Aut_C(\ssC,\ssL,\phi)$ and $\Aut'(C)$ (see also  \cite[chap. 2]{chiofar12}).
They fit into the following short exact sequence,
\begin{equation}\label{seq_aut}
1\to \underline\Aut_C(\ssC,\ssL,\phi)\to\underline\Aut(\ssC,\ssL,\phi)\to \Aut'(C)\to 1.
\end{equation}

\begin{defin}
Within a stable curve $C$, an elliptic tail is an irreducible component of geometric genus $1$ that
meets the rest of the curve in only one point called an elliptic tail node. Equivalently, $T$ is an elliptic tail if and only
if its algebraic genus is $1$ and $T\cap\overline{C\backslash T}=\{\msf{n}\}.$
\end{defin}
\begin{defin}\label{def_eti}
An element $i\in\Aut(C)$ is an elliptic tail automorphism if there exists
an elliptic tail $T$ of $C$ such that $i$ fixes $T$ and his restriction to $\overline{C\backslash T}$ is the identity.
An elliptic tail automorphism of order $2$ is called an elliptic tail quasireflection (ETQR). In the literature ETQRs are
called elliptic tail involutions (or ETIs), we changed this convention in order to generalize
the notion.
\end{defin}
\begin{rmk}\label{rmk_eti}
Every curve of algebraic genus $1$ with one marked point has
exactly one involution $i$. Then there is a unique ETQR associated to every elliptic tail.

More precisely an elliptic tail $E$ could be
of two types. The first type is a smooth curve of geometric genus $1$
with one marked point, \cor{i.e.} an elliptic curve: in this case we have $E=\C\slash\Lambda$,
for $\Lambda$ integral lattice of rank $2$, the marked point is the origin,
and the only involution is the map induced by $x\mapsto -x$ on $\C$.
The second type is the rational line with one marked point and an autointersection
point: in this case we can write $E=\p^1\slash\{1\equiv -1\}$, the marked point is the origin,
and the only involution is the map induced by $x\mapsto -x$ on $\p^1$.
\end{rmk}

\begin{teo}[{See \cite[theorem 2]{harmum82}}] \label{teo_eti}
Consider a stable curve $C$ of genus $g\geq 4$.
An element of $\Aut(C)$
acts as a quasireflection on $\defo(C)$ if and only if it is an ETQR. In particular, if $\eta\in \Aut(C)$ is an $ETQR$
acting non trivially on the tail $T$ with elliptic tail node $\msf{n}$, then 
 $\eta$ acts on $\defo(C)$ as 
$t_{\msf n}\mapsto -t_{\msf n}$ on the coordinate associated to $\msf{n}$, and
as $t\mapsto t$ on the remaining coordinates.
\end{teo}

\begin{rmk}\label{rmk_stacktail}
Consider a (stack-theoretic) curve $\msf E$ of genus $1$
with one marked point. We call $E$ its coarse space. In the case of an elliptic
tail of a curve $\ssC$, the marked point is the point of intersection between $\msf E$
and $\overline{\ssC\backslash \msf E}$.

If $\msf E$ is an elliptic curve, then $\msf E=E$ and
the curve has exactly one involution $i_0$.
In case $\msf E$ is rational, its normalization
is the stack $\overline{\msf E}=[\p^1\slash\mmu_r]$, with $\mmu_r$ acting by multiplication,
and $\msf E=\overline{\msf E}\slash\{0\equiv \infty\}$. There exists a canonical involution
$\msf i_0$ in this case too: we consider the pushforward of the involution
on $\p^1$ such that $z\mapsto 1\slash z$.

Given any twisted curve $\ssC$ with an elliptic tail $\msf E$ whose
elliptic tail node is $\msf n$,
we found a canonical involution $\msf i_0$ on $\msf E$
up to non-trivial action on $\msf n$.
\end{rmk}

\begin{defin}\label{def_etig}
We generalize the notion of ETQR to rooted curves.
An element $\msf i\in\Aut(\ssC,\ssL,\phi)$ is an ETQR if there exists
an elliptic tail $\msf E$ of $\ssC$ with elliptic tail node $\msf n$, such that
the action of $\msf i$ on $\ssC\backslash\msf E$ is trivial, and the action
on $\msf E$, up to non-trivial action on $\msf n$, is the canonical involution~$\msf i_0$.
\end{defin}

\section{The automorphisms of rooted curves}

In this section we will characterize smooth and singular points
of $\overline\R_{g,\ell}^k$ using the dual graph of curves and their
multiplicity indices. For simplicity, we will focus on the case $\ell$ prime.
See Remark \ref{rmk_composite} for a discussion on the generalization to composite level $\ell$.

Given a rooted curve $(\ssC,\ssL,\phi)$,
we come back to the local multiplicity on a node $\msf n$ whose
local picture is $[\{xy=0\}\slash\mmu_r]$ with $r$ positive integer.
As in equation \eqref{node_twist}, once we choose a privileged branch,
the action on the bundle fiber near
the node is $\xi_r(t,x,y)=(\xi_r^mt,\xi_r x,\xi_r^{-1}y)$.
We observe that the canonical line bundle $\omega_\ssC$ is the pullback
of the canonical line bundle over the coarse space $C$, and
this, with the isomorphism $\ssL^{\xx \ell}\cong \omega_\ssC$,
implies that ${(\xi_r^m)}^\ell=1$. So $\ell m$ is a multiple of~$r$.
As a consequence of the faithfulness of $\ssL$, the order $r$ 
equals $1$ or $\ell$.

\begin{defin}\label{mul_index}
The multiplicity index of $(\ssC,\ssL,\phi)$ is the cochain $M\in C^1(\Gamma;\Z\slash\ell)$
such that, for all $e\in\E$ oriented edge of the dual graph $\Gamma(\ssC)$,
$M(e)=m$ the local multiplicity
with respect
to the privileged branch associated to $e$ (see Remark \ref{rmk_or}).
\end{defin}
\begin{rmk}
If the node has trivial stabilizer, \cor{i.e.}~$r=1$, we pose $M(e):=0$.
\end{rmk}

We show a restating of Remark \ref{prop_count}.
\begin{defin}
Consider a
curve $C$ with dual graph $\Gamma$ and line bundle $F$ on $C$.
The multidegree cochains 
 is the function $\mdeg(F,C)\in C^0(\Gamma;\Z\slash\ell)$ such that
$$\mdeg(F,C)(v):=\deg F|_v\mod\ell \ \ \forall v\in V(\Gamma)$$
where by $F|_v$ we mean the restriction of $F$ on the
component associated to vertex $v$. 
\end{defin}
\begin{prop}\label{prop_mdeg}
Consider a stable curve $C$ with dual graph $\Gamma$ and consider a  cochain
$M$ in $C^1(\Gamma;\Z\slash\ell)$. Also consider the differential
$\partial\colon C^1(\Gamma;\Z\slash\ell)\to C^0(\Gamma;\Z\slash\ell)$
(see Section \ref{coch_sect}). There exists an $\ell$th root of the line bundle
$F$ on $C$ with multiplicity index $M$, if and only if 
$$\partial M\equiv\mdeg(F,C)\mod\ell.$$
\end{prop}
\proof If $(\ssC,\ssL,\phi)$ is any curve with an $\ell$th root of $F$, then we obtain the result
simply verifying a degree condition on every irreducible component of $\ssC$.
Consider the setting introduced in Remark~\ref{defo}: we call $v_i$ the $\Gamma$ vertex
associated to the irreducible components $\ssC_i$ and $g_i$ the genus of the normalization
$\overline\ssC_i$. Then we have
$$\deg\ssL|_{v_i}=r_i+\sum_{e_+=v_i}\frac{M(e)}{\ell},$$
with $r_i$ an integer for all $i$.
Knowing that $\left(\ssL|_{v_i}\right)^{\xx\ell}\cong F|_{v_i}$ and multiplying by $\ell$ 
the previous equality, we obtain, as we wanted,
$$\deg F|_{v_i}\equiv \sum_{e_+=v_i}M(e)\mod\ell.$$

To prove the other implication we will show that the
multidegree condition implies  $\deg F\equiv 0\mod\ell$, and
conclude by Remark \ref{prop_count}. Indeed,
$\deg F=\sum_{v_i\in V(\Gamma)}\deg F|_{v_i}$. Using
the condition we obtain
$$\deg F\equiv \sum_{v_i\in V(\Gamma)}\sum_{e_+=v_i}M(e)\equiv\sum_{e\in \E(\Gamma)}M(e)\equiv 0\mod\ell.$$\fine

Consider a rooted curve $(\ssC,\ssL,\phi)$ such that the coarse space
of $\ssC$ is $C$.
Starting from the dual graph $\Gamma(\ssC)$ and the multiplicity index $M$ of $(\ssC,\ssL,\phi)$,
consider the new contracted
graph $\Gamma_0(\ssC)$ defined~by
\begin{enumerate}
\item the vertex set $V_0=V(C)\slash\sim$, defined by modding out the relation $$(e_+\sim e_-\m{ if }M(e)\equiv0);$$
\item the edge set $E_0=\{e\in E(C)|\ M(e)\not\equiv 0\}$.
\end{enumerate}
\begin{rmk}
The graph $\Gamma_0$ is obtained by contracting the edges of $\Gamma$
where the function $M$ vanishes.
\end{rmk}
\begin{defin}\label{contr_graph}
The pair $(\Gamma_0(\ssC),M)$, where $M$ is the restriction
of the multiplicity index on the contracted edge set, is called
\emph{decorated graph} of the curve $(\ssC,\ssL,\phi)$. If the
cochain $M$ is clear from context, we will refer also to $\Gamma_0(\ssC)$ or $\Gamma_0$ alone
as the decorated graph.
\end{defin}

To study $\underline\Aut_C(\ssC,\ssL,\phi)$ we start
from a bigger group, the group $\Aut_C(\ssC)$ containing automorphisms of $\ssC$ fixing the
coarse space $C$.
Consider   a node $\msf{n}$ of $\ssC$ whose local picture is $[\{xy=0\}\slash\mmu_r]$.
Consider an automorphism $\eta\in\Aut_C(\ssC)$.
The local action of $\eta$ at $\msf{n}$ is
$(x,y)\mapsto(\zeta x,y)=(x,\zeta y)$,
with $\zeta\in\mmu_r$.
As a consequence of the definition of $\Aut_C(\ssC)$, the action
of $\eta$ outside the $\ssC$ nodes is trivial. 
Then the whole
group $\Aut_C(\ssC)$ is generated by automorphisms of the form
$(x,y)\mapsto (\zeta x, y)$ on a node and trivial elsewhere.

We are interested in representing $\Aut_C(\ssC)$ as acting on the edges of the dual graph, thus
we introduce the group of functions $\E\to\Z\slash\ell$ that are even with respect to conjugation
$$S(\Gamma;\Z\slash\ell):=\left\{b\colon\E\to \Z\slash\ell\ |\ b(\bar e)\equiv b(e)\right\}.$$
We have a canonical identification sending the function $b\in S(\Gamma_0(\ssC);\Z\slash\ell)$
to the automorphism $\eta$ with local action $(x,y)\mapsto(\xi_\ell^{b(e)} x,y)$ on the node 
associated to the edge $e$ if $M(e)\not\equiv 0$. Therefore
the decorated graph encodes the automorphisms acting
trivially on the coarse space. We can write
\begin{equation}\label{aut=Snu}
\Aut_C(\ssC)\cong S(\Gamma_0(\ssC);\Z\slash\ell)\cong (\Z\slash\ell)^{E_0(\Gamma)}.
\end{equation}
We already saw that elements of $C^1(\Gamma;\Z\slash\ell)$ are odd functions
from $\E$ to $\Z\slash\ell$. Given $g\in S(\Gamma;\Z\slash\ell)$ and $N\in C^1(\Gamma;\Z\slash\ell)$,
their natural product $gN$ is still an odd function, thus an element of $C^1(\Gamma;\Z\slash\ell)$.

Given the normalization morphism
$\msf{nor}\colon C_{\msf{nor}}\to C$, consider 
the short exact sequence of sheaves over $C$
$$ 1\to \Z\slash\ell\to {\msf {nor}}_*{\msf {nor}}^*\Z\slash\ell\xrightarrow{ t } \left.\Z\slash\ell\right|_{\sing C}\to 1.$$
The sections of the central sheaf are $\Z\slash\ell$-valued functions over $C_{\msf{nor}}$.
Moreover, 
the image of a section $s$ by $t$ is the function
that assigns to each node the difference between the two values
of $s$ on the preimages.
The cohomology of this sequence gives the following
long exact sequence
\begin{equation}\label{es_tau}
1\to \Z\slash\ell\to C^0(\Gamma;\Z\slash\ell)\xrightarrow{\ \delta\ }C^1(\Gamma;\Z\slash\ell)\xrightarrow{\ \tau\ }\Pic(C)[\ell]
\xrightarrow{\msf{nor}^*}\Pic(C_{\msf{nor}})[\ell]\to 1.
\end{equation}
Here, $\Pic(C)[\ell]=H^1(C,\Z\slash\ell)$ is the subgroup of $\Pic(C)$ of elements of order dividing $\ell$, \cor{i.e.}~of 
$\ell$th roots of the trivial bundle.

We know that
$$\underline\Aut_C(\ssC,\ssL,\phi)=\{[\msf{s},\rho]\in\underline\Aut(\ssC,\ssL,\phi)|\ \msf{s}\in\Aut_C(\ssC)\}=\{\msf{s}\in\Aut_C(\ssC)|\ \msf{s}^*\ssL\cong \ssL\}.$$
As showed in \cite{chio08}, we have
$$\msf{a}^*\ssL=\ssL\xx\tau\left(\msf{a}M\right).$$
This means that,
if $\msf{a}$ is an automorphism in $\Aut_C(\ssC)$,
the pullback $\msf{a}^*\ssL$ is totally determined by the product $\msf{a}M\in\C^1(\Gamma_0;\Z\slash\ell)$,
where $\msf{a}$ is seen as
an element of $S(\Gamma_0;\Z\slash\ell)$, and
 $M$ is the multiplicity index of $(\ssC,\ssL,\phi)$.

As a consequence we have the following theorem.
\begin{teo}\label{lift_condition}
An element $\msf{a}\in\Aut_C(\ssC)$, lifts to $\underline\Aut_C(\ssC,\ssL,\phi)$
if and only if
$$\msf{a}M\in\Ker(\tau)=\im\delta.$$
\end{teo}
We recall the subcomplex $C^i(\Gamma_0(\ssC);\Z\slash\ell)\subset C^i(\Gamma(\ssC);\Z\slash\ell)$.
Moreover, if $\delta_0$ is the $\delta$ operator on $C^i(\Gamma_0)$,
\cor{i.e.}~the restriction of the $\delta$ operator to this space, 
 from Proposition \ref{prop_imcap} we know that
$\im(\delta_0)=C^1(\Gamma_0)\cap\im\delta.$
\begin{rmk}\label{rmk_G}
As a consequence
of Theorem \ref{lift_condition}, we get the canonical identification
$$\underline\Aut_C(\ssC,\ssL,\phi) = \im(\delta_0)$$ 
via the multiplication $\msf a\mapsto \msf aM$. Because of Proposition \ref{prop_G} we also
have $\underline\Aut_C(\ssC,\ssL,\phi)=M\cdot G(\Gamma_0;\Z\slash\ell)$.
\end{rmk}
\begin{rmk}\label{rmk_prod}
In the first section we obtained a characterization of the cochains
in $\im(\delta)$ that we could restate in our new setting. Indeed,
because of Proposition \ref{prop_imdelta}, an automorphism $\msf{a}\in S(\Gamma_0;\Z\slash\ell)$ is
an element of $\underline\Aut_C(\ssC,\ssL,\phi)$ if and only if for every circuit $(e_1,\dots,e_k)$ in $\Gamma_0$
we have $\sum_{i=1}^k \msf a(e)\equiv 0\mod \ell$.
\end{rmk}

We can characterize the automorphisms in $\underline\Aut_C(\ssC,\ssL,\phi)$
that are quasireflections.
We will regard every automorphism $\msf{a}\in\underline\Aut_C(\ssC,\ssL,\phi)$ as
the corresponding element in $S(\Gamma_0;\Z\slash\ell)$. This representation 
is good to investigate the action of $\msf{a}$ on $\defo(\ssC)$. The action
of  $\msf{a}$ is
non-trivial only on those coordinates of $\defo(\ssC)$ associated to the nodes of $\ssC$ with non-trivial stabilizer. In particular,
if $t_{\msf{n}}$ is the coordinate associated to the node $\msf{n}$ and $e$ is the edge
of $\Gamma_0(\ssC)$ associated to $\msf{n}$, then 
$$\msf{a}\colon t_{\msf{n}}\mapsto \xi_\ell^{\msf{a}(e)}\cdot t_{\msf{n}}.$$

\begin{prop}\label{ghost1}
An automorphism $\msf{a}\in\Aut_C(\ssC)$ is a quasireflection in $\underline\Aut_C(\ssC,\ssL,\phi)$ if and only if
$\msf{a}(e)\equiv 0$ for all edges but one  that is a separating edge of $\Gamma_0(\ssC)$.
\end{prop}

\proof If $\msf{a}$ is a quasireflection, the action on all but one of the coordinates must be trivial.
Therefore $\msf{a}(e)\equiv 0$ on all the edges but one, say $e_1$.
If $e_1$ is in any circuit $(e_1,\dots,e_k)$ of $\Gamma_0$ with $k\geq 1$, we have, by Remark \ref{rmk_prod},
that $\sum \msf{a}(e_i)\equiv 0$.
As $\msf{a}(e_1)\not\equiv 0$, there exists $i>1$ such that $\msf{a}(e_i)\not\equiv 0$, contradiction.
Thus $e_1$ is not in any circuit, then it is a separating edge. 

Reciprocally,
consider an automorphism $\msf{a}\in S(\Gamma_0;\Z\slash\ell)$ such
that there exists an oriented separating edge $e_1$ with the property that
$\msf{a}(e)\equiv 0$ for all $e\in \E\backslash\{e_1, \bar e_1\}$ and
$\msf{a}(e_1)\not\equiv 0$. Then for every circuit $(e_1',\dots,e_k')$ we have
$\sum \msf{a}(e_i')\equiv 0$ and so $\msf{a}$ is in $\underline\Aut_C(\ssC,\ssL,\phi)$.\fine

\begin{defin}
We call $\QR(\underline\Aut(\ssC,\ssL,\phi))$ or simply $\QR(\ssC,\ssL,\phi)$
the group generated by quasireflections automorphism of the curve $(\ssC,\ssL,\phi)$.
We call $\QR(\underline\Aut_C(\ssC,\ssL,\phi))$ or simply $\QR_C(\ssC,\ssL,\phi)$
the group generated by ghost quasireflections.
\end{defin}

\begin{rmk}
If we note $E_{\sep}\subset E_0$ the subset of separating edges of $\Gamma_0$, 
Proposition \ref{ghost1} gives a simple description of the group $\QR_C(\ssC,\ssL,\phi)$.
$$\QR_C(\ssC,\ssL,\phi)\cong(\Z\slash\ell)^{E_{\sep}}\subset S(\Gamma_0 ;\Z\slash\ell).$$
\end{rmk}

\begin{teo}\label{treelike1}
The group $\underline\Aut_C(\ssC,\ssL,\phi)$ is generated by quasireflections
if and only if the graph $\Gamma_0$ is tree-like.
\end{teo}
\proof After the previous remark,
$$\QR_C(\ssC,\ssL,\phi)\cong(\Z\slash\ell)^{E_{\sep}}\subset\im(\delta_0)\cong\underline\Aut_C(\ssC,\ssL,\phi).$$
We know that $\im(\delta_0)\cong\Z\slash\ell^{\#V_0-1}$ and thus the inclusion
is an equality if and only if $\#E_{\sep}=\#V_0-1$. By Lemma \ref{ineq_graph2}
we conclude.\fine

\begin{rmk}\label{rmk_composite}
We show a generalization of the previous result to every $\ell$,
and to do this we update our tools following \cite[\S2.4]{chiofar12}.
First of all we generalize the notion of multiplicity index:
given the local multiplicity $m\in\Z\slash r$ at a certain node $\msf n$,
we know that $r$ divides $m\ell$, we define $M(e):=m\ell\slash r\in \Z\slash\ell$,
where $e\in\E$ is an oriented edge. If $\ell$ is prime, 
$M\in C^1(\Gamma;\Z\slash\ell)$ coincide
with the previous definition.
We also generalize contracted graphs.
Consider the factorization of $\ell$ in prime numbers,
$\ell=\prod p^{e_p}$, where $e_p=\nu_p(\ell)$ is the $p$-adic
valuation of $\ell$.
Given the dual graph $\Gamma(\ssC)=\Gamma$,
for every prime $p$
we define $\Gamma(\nu_p^k)$ contracting every edge $e$ of
$\Gamma$ such that $p^k$ divides $M(e)$. 
We will note $\Gamma_p(\ssC):=\Gamma(\nu_p^{e_p})$.
We have the chain of contraction already introduced in \cite{chiofar12}:
$$\Gamma\to\Gamma_p=\Gamma(\nu_p^{e_p})\to\Gamma(\nu_p^{e_p-1})\to\cdots\to\Gamma(\nu_p^1)\to\{\cdot\}.$$
Moreover, we introduce $S_d(\Gamma;\Z\slash\ell)\subset S(\Gamma;\Z\slash\ell)$, the 
group of even functions $f\colon\E\to\Z\slash\ell$ such that $f(e)=f(\bar e)\in\Z\slash r(e)$,
and $C^1_d(\Gamma;\Z\slash\ell)\subset C^1(\Gamma;\Z\slash\ell)$, the group of odd functions $N\colon\E\to\Z\slash\ell$
such that $N(e)=-N(\bar e)\in\Z\slash r(e)$. We point out that if $r$ divides $\ell$, we will
look at $\Z\slash r$ as immersed in $\Z\slash\ell$, indeed
$\Z\slash r=(\ell\slash r)\cdot\Z\slash r\subset\Z\slash\ell$.
In particular for every prime~$p$,
$\Z\slash p\subset\Z\slash p^2\subset\Z\slash p^3\subset\cdots$.

\begin{teo}\label{teo_composite}
Consider a rooted curve $(\ssC,\ssL,\phi)$. 
The groups $\underline\Aut_C(\ssC,\ssL,\phi)$ is generated by quasireflections
if and only if the graphs $\Gamma_p(\ssC)$ are tree-like for every prime
$p$ dividing~$\ell$.
\end{teo}
\proof As in the case with $\ell$ prime,
$S_d(\Gamma;\Z\slash\ell)$ is canonically identified with $\Aut_C(\ssC)=\bigoplus_{e\in E}\Z\slash r(e)$.
Given the multiplicity index $M$ of $(\ssC,\ssL,\phi)$,
as before we have that $\msf a$ lifts to $\underline\Aut_C(\ssC,\ssL,\phi)$
if and only if $\msf aM$ is in $\im\delta$. Thus $\underline\Aut_C(\ssC,\ssL,\phi)$
is canonically identified with $C_d^1(\Gamma;\Z\slash\ell)\cap\im\delta$.

We introduce the $p$-adic valuation on edges as in \cite[\S 2.4.2]{chiofar12}:
given $e\in\E$, $\nu_p(e):=\m{val}_p(M(e)\mod p^{e_p})$.
We consider also the function $\bar \nu_p$ such that
 $\bar \nu_p(e):=\min(e_p,\nu_p(e))$. We observe that, by definition
$M$, for every oriented edge $e\in\E(\Gamma)$ the local stabilizer order is
$$r(e)=\prod_{p|\ell}p^{e_p-\bar\nu_p(e)}.$$
We define one last new object. We have a canonical immersion
$C^i(\Gamma(\nu_p^k);\Z\slash p^{e_p-k+1})\hookrightarrow C^i(\Gamma;\Z\slash\ell)$
for $i=0,1$. We define $\delta_p^k$ as the restriction of the delta operator on the
group $C^0(\Gamma(\nu_p^k);\Z\slash p^{e_p-k+1})$
for all $p$ in the factorization of $\ell$ and for all $k$ between $1$ and $e_p$.

As in Lemma \ref{ghost1}, in the composite case $\msf a\in\Aut_C(\ssC)$
is a quasireflection in $\underline\Aut_C(\ssC,\ssL,\phi)$ if and only if $\msf a (e)\equiv 0$
for all edges but one that is a separating edge. This allows the following decomposition
$$\QR_C(\ssC,\ssL,\phi)=\bigoplus_{e\in E_{\sep}(\Gamma)}\Z\slash r(e)=\bigoplus_{e\in E_{\sep}}
\bigoplus_{p|\ell}\Z\slash p^{e_p-\bar\nu_p(M(e))}=\bigoplus_{p|\ell}\bigoplus_{k=1}^{e_p}
(\Z\slash p^k)^{\beta_p^k},$$
where $\beta_p^k:=\# E_{\sep}\left(\Gamma(\nu_p^{e_p-k+1})\right)-\# E_{\sep}\left(\Gamma(\nu_p^{e_p-k})\right)$ if $k<e_p$
and $\beta_p^{e_p}:=\# E_{\sep}\left(\nu_p^1\right)$.
Following \cite[Lemma 2.22]{chiofar12}, we have a similar decomposition for
$\underline\Aut_C(\ssC,\ssL,\phi)$:
$$\underline\Aut_C(\ssC,\ssL,\phi)=C^1_d(\Gamma;\Z\slash\ell)\cap\im\delta=\bigoplus_{p|\ell}\sum_{k=1}^{e_p}\im\delta_p^k=\bigoplus_{p|\ell}
\bigoplus_{k=1}^{e_p}(\Z\slash p^k)^{\alpha_p^k},$$
where $\alpha_p^k:=\# V\left(\Gamma(\nu_p^{e_p-k+1})\right)-\# V\left(\Gamma(\nu_p^{e_p-k})\right)\ \forall k\geq 0.$
We observe that $\alpha_p^k\geq \beta_p^k$ for all $p$ dividing $\ell$ and $k\geq 0$. Moreover,
$\underline\Aut_C(\ssC,\ssL,\phi)$ coincide with $\QR_C(\ssC,\ssL,\phi)$
 if and only if $\alpha_p^k=\beta_p^k$ for all $p$ and $k$. Fixing
$p$, this is equivalent to impose $\sum_k\beta_p^k=\sum_k\alpha_p^k$. In the previous expression
the left hand side is $\#E_{\sep}(\Gamma_p)$ and the right hand side is $\#V(\Gamma_p)-1$,
we saw in Lemma \ref{ineq_graph2} that the equality is achieved if and only if $\Gamma_p$ is tree-like.\fine
\end{rmk}

\section{The singular locus of $\overline{\R}_{g,\ell}^k$}

After the analysis of quasireflections on $\underline\Aut_C(\ssC,\ssL,\phi)$,
we complete the description of quasireflections on the whole automorphism groups. What
follows is true for every $\ell$.

\begin{lemma}\label{lemma_qr}
Consider an element $\msf q$ of $\underline\Aut(\ssC,\ssL,\phi)$. It acts as a quasireflection on $\defo(\ssC)$
if and only if one of the following is true:
\begin{enumerate}
\item the automorphism $\msf q$ is a ghost quasireflection,
\cor{i.e.}~an element of $\underline\Aut_C(\ssC,\ssL,\phi)$ which moreover operates as
a quasireflection;
\item  the automorphism $\msf q$ is an ETQR, using the generalized Definition \ref{def_etig}.
\end{enumerate}
\end{lemma}
\proof 
We first prove the ``only if'' part. Consider $\msf q\in\underline\Aut(\ssC,\ssL,\phi)$, we call $q$ its coarsening.
If $\msf q$ acts trivially on certain coordinates of $\defo(\ssC)$, \cor{a fortiori} we have that $q$
acts trivially on the corresponding coordinates of $\defo(C)$. So $q$ acts as the identity or as a quasireflection.
In the first case, $\msf q$ is a ghost automorphism and we are in case $(1)$. 
If $q$ acts as
a quasireflection, then it is a classical ETQR as we pointed out on Theorem \ref{teo_eti},
and it acts non-trivially on the elliptic tail node $\msf n$ associated to an elliptic tail.
It remains to know the action of $\msf q$ on the nodes, other than $\msf n$, with non-trivial stabilizer.
The action in these nodes must be trivial, because every non-trivial action induces a non-trivial action
also on the associated
coordinate of the node. Therefore, the $\msf q$ restriction to the elliptic tail
has to be the canonical involution $\msf i_0$ (see Remark \ref{rmk_stacktail}).
By Definition \ref{def_etig} this implies that $\msf q$ is an ETQR of $(\ssC,\ssL,\phi)$.

For the ``if'' part, we observe that a ghost quasireflection is automatically a quasireflection.
It remains to prove the point $(2)$. By definition of ETQR its action can be non-trivial
only on the separating node of the tail. The local coarse picture of the node is $\{xy=0\}$,
where $y=0$ is the branch lying on the elliptic tail. Then the action of $\msf i$ on the coarse space
is $(x,y)\mapsto(-x,y)$. Therefore the action is \cor{a fortiori} non trivial
on the coordinate associated to the stack node $\msf n$.\fine

\begin{lemma}\label{keylemma}
If $\QR(\Aut(C))$ (also called $\QR(C)$) is the group generated
by ETQRs inside $\Aut(C)$, then any element $q\in \QR(C)$ which
could be lifted to $\underline\Aut(\ssC,\ssL,\phi)$, has a lifting in~$\QR(\ssC,\ssL,\phi)$, too.
\end{lemma}
\proof By definition, $\underline\Aut(\ssC,\ssL,\phi)$ is the set of automorphisms $\msf s\in\Aut(\ssC)$
such that $\msf s^*\ssL\cong\ssL$. Consider $q\in\QR(C)$ such that its decomposition
in quasireflections is $q=i_0i_1\cdots i_m$, and $i_k$ is an ETQR acting non-trivially on an elliptic tail $E_k$.
Any lifting of $q$ is in the form $\msf q=\msf i_0\msf i_1\cdots \msf i_m\cdot \msf a$, where 
$\msf i_k$ is a (generalized) ETQR acting non-trivially on $E_k$, and $\msf a$ is a ghost
acting non-trivially only on nodes outside the tails $E_k$. We observe that every $\msf i_k$ is a lifting
in $\Aut(\ssC)$ of $i_k$, we are going to prove that moreover $\msf i_k\in \underline\Aut(\ssC,\ssL,\phi)$.
By construction, $\msf q^*\ssL\cong\ssL$ if and only if $\msf i_k^*\ssL\cong \ssL$ for
all $k$ and $\msf a^*\ssL\cong \ssL$. This implies that every $\msf i_k$ lives
in $\underline\Aut(\ssC,\ssL,\phi)$, then $\msf q\msf a^{-1}$ is a lifting of $q$
living in $\QR(\ssC,\ssL,\phi)$.\fine

\begin{rmk}\label{keyrmk}
We recall the short exact sequence
$$1\to\underline\Aut_C(\ssC,\ssL,\phi)\xrightarrow{\ h\ }\underline\Aut(\ssC,\ssL,\phi)\xrightarrow{\ p\ }\Aut'(C)\to 1$$
and introduce the group $\QR'(C)\subset\Aut'(C)$, generated by liftable quasireflections.
Knowing that $p(\QR(\ssC,\ssL,\phi))\subset \QR'(C)\subset \Aut'(C)\cap\QR(C)$,
the previous lemma shows that $\QR'(C)=p(\QR(\ssC,\ssL,\phi))$. Using also Lemma \ref{lemma_qr},
we obtain that the following is a short exact sequence
$$1\to \QR_C(\ssC,\ssL,\phi)\to \QR(\ssC,\ssL,\phi)\to \QR'(C)\to 1.$$
\end{rmk}

\begin{teo}\label{teoqr}
The group $\underline\Aut(\ssC,\ssL,\phi)$ is generated by quasireflections
if and only if both $\underline\Aut_C(\ssC,\ssL,\phi)$ and $\Aut'(C)$ are generated by
quasireflections.
\end{teo}
\proof
After the previous remark, the following is a short exact sequence
$$1\to \underline\Aut_C(\ssC,\ssL,\phi)\slash\QR_C(\ssC,\ssL,\phi)\to
\underline\Aut(\ssC,\ssL,\phi)\slash\QR(\ssC,\ssL,\phi)\to \Aut'(C)\slash\QR'(C)\to 1.$$
The theorem follows.\fine

After Theorem \ref{teoqr} and \ref{treelike1}, we have the following characterization
of the singular locus inside $\overline\R_{g,\ell}^k$.

\begin{teo}\label{smooth2}
Let $g\geq 4$ and $\ell$ a prime number. Given a rooted curve $(\ssC,\ssL,\phi)$, with
$C$ coarse space of $\ssC$, the point $[\ssC,\ssL,\phi]\in\overline\R_{g,\ell}^k$ is smooth
if and only if $\Aut'(C)$ is generated by ETQRs of $C$ and the contracted graphs $\Gamma_0(\ssC)$
is tree-like.
\end{teo}

After Remark \ref{rmk_composite}, we can generalize the previous theorem to all $\ell$.
 It suffices
to consider contracted graphs $\Gamma_p(\ssC)$ for every prime $p$ dividing $\ell$
\begin{teo}\label{smooth1}
For any $g \geq 4$ and $\ell$ positive integer. The point $[\ssC,\ssL,\phi]$
is smooth if and only if $\Aut'(C)$ is generated by ETQRs of $C$ and the $\Gamma_p(\ssC)$ are tree-like.
\end{teo}
\begin{rmk}\label{rmk_teo}
The theorem above is a generalization of \cite[Theorem 2.28]{chiofar12}. In particular Chiodo and Farkas proved
that in the case $k=0$, the moduli of level curves, the contracted graphs $\Gamma_p$ are bouquet
if $[\ssC,\ssL,\phi]$ is smooth. A bouquet is a graph with only one vertex, or equivalently
a tree-like graph with no non-loop edges. This happens because in the case $k=0$, every separating
node must have trivial stabilizer, \cor{i.e.}~the $M$ cochain must be $0$ and therefore they disappear
after the contraction.
In the general case, $M$ could be non-zero on separating edges too, and the theorem
is still true using the more general notion of tree-like graph.
\end{rmk}

\subsection{The locus of singular points via a new stratification}\label{subs_strata}
As we saw, the information about the automorphism group of a certain rooted curve $(\ssC,\ssL,\phi)$
is coded by its dual decorated graph $(\Gamma_0(\ssC),M)$. It is therefore quite natural to introduce
a stratification of $\overline\R_{g,\ell}^k$ using this notion.
For this, we extend the notion of
graph contraction: if $\Gamma_0'\to\Gamma_1'$ is a usual graph contraction,
the ring $C^1(\Gamma_1';\Z\slash\ell)$ is naturally immersed in $C^1(\Gamma_0';\Z\slash\ell)$,
then the contraction of a pair $(\Gamma_0',M_0')$ is the pair $(\Gamma_1',M_1')$ where
the cochain $M_1'$ is the restriction of $M_0'$. If it is clear from the context,
we could refer to the decorated graph restriction simply
with the graph contraction $\Gamma_0'\to\Gamma_1'$.

\begin{defin}
Given a decorated graph $(\Gamma,M)$ with $M\in C^1(\Gamma;\Z\slash\ell)$, consider this
locus of $\overline\R_{g,\ell}^k$:
$$\s_{(\Gamma,M)}:=\left\{[\ssC,\ssL,\phi]\in\overline\R_{g,\ell}^k:\ \Gamma_0(\ssC)=\Gamma,\ \m{and }
M\ \m{is the multiplicity index of }(\ssC,\ssL,\phi)\right\}.$$
\end{defin}
These loci are the strata of our stratification. We can find a first
link between the decorated graph and geometric properties of the stratum.
\begin{prop}\label{prop_codim}
If we consider the codimension of $\s_{(\Gamma,M)}$ inside $\overline\R_{g,\ell}^k$, we have
$$\codim \s_{(\Gamma,M)}=\#E(\Gamma).$$
\end{prop}
\proof We take a general point $[\ssC,\ssL,\phi]$ of stratum $\s_{(\Gamma,M)}$, it has
$\#V(\Gamma)$ irreducible components $\ssC_1,\dots,\ssC_{\#V}$. We call
$g_i$ the genus of $\ssC_i$ and $k_i$ the number of nodes of $\ssC$ on this component. Then
we have $\sum k_i=2\#E(\Gamma)$.
We obtain that the dimension of $\defo(\ssC,\ssL,\phi)$ is
$$\dim\defo(\ssC,\ssL,\phi)=\sum_{i=1}^{\#V}(3g_i-3+k_i)=3\sum g_i-3\#V+2\#E.$$
Thus, because of Remark \ref{gen_nodal}, we have 
$\dim\defo(\ssC,\ssL,\phi)=3g-3-\#E$, where $g$ is the genus of $\ssC$. The result
on the codimension follows.\fine

Using contraction, we have this
description of the closed strata.
$$\overline\s_{(\Gamma_1',M_1')}=\left\{[\ssC,\ssL,\phi]\in\overline\rr_{g,\ell}^k:\ 
\begin{minipage}{8cm}
if $(\Gamma_0, M)$ is the decorated graph of 
$(\ssC,\ssL,\phi)$, there exists a contraction $(\Gamma_0,M)\to(\Gamma_1',M_1')$
\end{minipage}
\right\}.$$

We introduce two closed loci of $\overline\R_{g,\ell}^k$,
$$N_{g,\ell}^k:=\left\{[\ssC,\ssL,\phi]|\ \Aut'(C)\m{ is not generated by ETQRs}\right\},$$
$$H_{g,\ell}^k:=\left\{[\ssC,\ssL,\phi]|\ \underline\Aut_C(\ssC,\ssL,\phi)\m{ is not generated by quasireflections}\right\}.$$
Equivalently, $[\ssC,\ssL,\phi]\in H_{g,\ell}^k$ if and only if  $\Gamma_0(\ssC)$ is not tree-like.
We have by Theorem \ref{smooth2} that the singular locus $\sing \overline\R_{g,\ell}^k$ is
their union
$$\sing\overline\R_{g,\ell}^k=N_{g,\ell}^k\cup H_{g,\ell}^k.$$

\begin{rmk}
Consider the natural projection  $\pi\colon\overline\R_{g,\ell}^k\to\overline\M_g$,
we observe that 
$$N_{g,\ell}^k\subset\pi^{-1}\sing\overline\M_g.$$
Indeed, after Remark \ref{keyrmk}, $\QR'(C)=\Aut'(C)\cap \QR(C)$
and therefore $\Aut(C)=\QR(C)$ if and only if $\Aut'(C)=\QR'(C)$.
This implies that $\left(\pi^{-1}\sing\overline\M_g\right)^c\subset (N_{g,\ell}^k)^c$,
and taking the complementary we have the result.
\end{rmk}

The stratification introduced is particularly useful in describing 
the ``new'' locus $H_{g,\ell}^k$. We recall the definition of \emph{vine} curves.

\begin{defin}\label{vine}
We note $\Gamma_{(2,n)}$ a graph with two vertices linked by $n$ edges. A curve $\ssC$
is an $n$-vine curve if $\Gamma(\ssC)$ contracts to $\Gamma_{(2,n)}$ for some $n\geq 2$.
Equivalently an $n$-vine curve has coarse space $C=\ssC_1\cup C_2$ which is
the union of two curves intersecting each other $n$ times. 
\end{defin}

The next step shows that every $H$-curve has
a dual graph which contracts to $\Gamma_{(2,n)}$.
\begin{lemma}
If $[\ssC,\ssL,\phi]\in\sing\overline\R_{g,\ell}^k$ is a point in $H_{g,\ell}^k$
whose decorated graph is $(\Gamma_0,M)$,
then there exists $n\geq 2$ and a graph contraction $\Gamma_0\to \Gamma_{(2,n)}$.
Equivalently $\ssC$ is an $n$-vine curve for some $n\geq 2$.
\end{lemma}
\proof 
If $[\ssC,\ssL,\phi]\in H_{g,\ell}^k$, by definition $\Gamma_0$ contains a circuit that is not a loop.
We will show an edge contraction $\Gamma_0\to \Gamma_{(2,n)}$ 
for some $n\geq 2$.

Consider two different vertices $v_1$ and $v_2$ that are consecutive on a non-loop circuit $K\subset\Gamma_0$.
Now consider a partition $V=V_1\bigsqcup V_2$ of the vertex set such that
$v_1	\in V_1$ and $v_2\in V_2$. This defines an edge contraction $\Gamma_0\to \Gamma_{(2,n)}$
where $e\in E(\Gamma)$ is contracted if its two extremities lie in the same $V_i$. As $K$ is
a circuit, necessarily $n\geq 2$ and the theorem is proved.\fine

We conclude that $H_{g,\ell}^k$ is the union of the closed strata associated
to vine graphs.
$$H_{g,\ell}^k\ \ =\bigcup_{
\substack{ n\geq 2\\
M\in C^1(\Gamma_{(2,n)};\Z\slash\ell)}}\overline\s_{(\Gamma_{(2,n)},M)}.$$

\section{Non-canonical singularities of $\overline\R_{g,\ell}^k$}

After the description of the singular locus of $\overline\R_{g,\ell}^k$, in this
section we find the locus of non-canonical singularities $\sing^{\nc}\overline\R_{g,\ell}^k$. To do this
we introduce the
 $\age$ function and point out some basic facts
about quotient singularities: the characterization of the non-canonical locus
will follow from the age criterion \ref{age_crit}.
Also, for some small values of $\ell$,
 we will develop a description of $\sing^{\nc}\overline\R_{g,\ell}^k$ in terms
of the stratification introduced above. 

\subsection{The age criterion}
If $G$ is a finite group, the age is an additive positive function from
the representations ring of $G$ to rational numbers,
$$\age\colon \Rep(G)\to \Q.$$
First consider the group $\Z\slash r$ for any $r$ positive integer. 
Given the character $\pmb k$ such that $1 \mapsto k\in\Z\slash r$,
we define $\age(\pmb{k})=k\slash r$. These
characters are a basis for $\Rep(\Z\slash r)$, then we can extend age
over all the representation ring.

Consider a $G$-representation $\rho\colon G\to \GL(V)$, the age function
could be defined on any injection $i\colon \Z\slash r\hookrightarrow G$ simply composing the
injection with $\rho$.
$$\age_V\colon i\mapsto \age(\rho\circ i).$$
Age can finally be defined on the group $G$ by
$$G\xrightarrow{\ f\ }\bigsqcup_{r\geq 1}\{i|\ i\colon \mmu_r\hookrightarrow G\}\xrightarrow{\age_V}\Q,$$
where $f$ is the set bijection sending $g\in G$, element of order $r$, to the injection
obtained by mapping $1\in\Z\slash r$ to $g$.
\begin{defin}[junior and senior groups]\label{def_junior}
A finite group $G\subset \GL(\C^m)$ that contains no quasireflections
 is called junior if the image of the age function intersects the open
interval~$]0,1[$,
$$\age G\cap\ ]0,1[\ \neq \varnothing,$$
and senior otherwise.
\end{defin}
\begin{rmk}\label{rmk_choice}
All the functions we have introduced are canonically defined except for the identification $f$,
that depends on the choice of the primitive root $\xi_r=\exp(2\pi i\slash r)$.
Moreover, the image of $\age\colon G\to \Q$ does not depend on this choice, then
our previous choice of a root system does not affect the study of junior and senior groups
and will continue unchanged.
\end{rmk}

\begin{teo}[age criterion, \cite{reid80}]\label{age_crit}
If $G\subset\GL(\C^m)$ is a finite subgroup without quasireflections,
the singularity $\C^m\slash G$ is canonical if and only if $G$ is a senior group.
\end{teo}

\begin{rmk}
To see the age explicitly, for $g\in G\subset\GL(\C^m)$ with $G$ finite subgroup and
$\ord(g)=r$, consider a basis of $\C^m$ such that $g=\diag(\xi_r^{a_1},\dots,\xi_r^{a_m})$. 
In this setting $\age(g)=\frac{1}{r}\sum a_i$.
\end{rmk}

\subsection{The non-canonical locus}

We know that any point $[\ssC,\ssL,\phi]\in\overline\R_{g,\ell}^k$ has a neighborhood
isomorphic to the quotient
 $\defo(\ssC)\slash\underline\Aut(\ssC,\ssL,\phi)$, then after the age criterion \ref{age_crit}
we are searching for junior automorphism in $\underline\Aut(\ssC,\ssL,\phi)$.
We also need the denominator group to be quasireflections free to use the criterion,
so we will repeatedly use the following result.
\begin{prop}[see \cite{prill67}]\label{propt}
Consider the finite subgroup $G\subset\GL(\C^m)$ and the group $\QR(G)$ generated
by $G$ quasireflections. There exists an isomorphism $\varphi\colon \C^m\slash \QR(G)\to\C^m$
and a subgroup $K\subset\GL(\C^m)$ isomorphic to $G\slash\QR(G)$ such that the following
diagram is commutative.
\[
\begin{CD}
\C^m@>>> \C^m\slash \QR(G)@>\varphi>>\C^m\\
@VVV @VVV @VVV\\
\C^m\slash G@>\cong>> (\C^m\slash \QR(G))\slash (G\slash\QR(G))@>\cong>> \C^m\slash K\\
\end{CD}
\]\newline
\end{prop}

We introduce two closed loci which are central
in our description.
\begin{defin}[$T$-curves]
A rooted curve $(\ssC,\ssL,\phi)$ is a $T$-curve if
there exists
an automorphism $\msf a\in \underline\Aut(\ssC,\ssL,\phi)$ such that its 
coarsening $a$ is an elliptic tail automorphism of order $6$.
The locus of $T$-curves in $\overline\R_{g,\ell}^k$ is noted $T_{g,\ell}^k$.
\end{defin}

\begin{defin}[$J$-curves]
A rooted curve $(\ssC,\ssL,\phi)$ is a $J$-curve if
the group 
$$\underline\Aut_C(\ssC,\ssL,\phi)\slash\QR_C(\ssC,\ssL,\phi),$$
which is the group of ghosts quotiented by its subgroup of quasireflections, is junior.
The locus of $J$-curve in $\overline\R_{g,\ell}^k$ is noted $J_{g,\ell}^k$.
\end{defin}

\begin{teo}\label{teo_jt}
For $g\geq 4$, the non-canonical locus of $\overline\R_{g,\ell}^k$ is formed
by $T$-curves and $J$-curves, \cor{i.e.}~it is the union
$$\sing^{\nc}\overline\R_{g,\ell}^k=T_{g,\ell}^k\cup J_{g,\ell}^k.$$
\end{teo}

\begin{rmk}
We observe that Theorem 2.44 of Chiodo and Farkas \cite{chiofar12}, affirms
exactly that in the case $k=0,\ \ell\leq 6$ and $ \ell \neq 5$, the $J$-locus
$J_{g,\ell}^0$ is empty for every genus $g$, and therefore $\sing^{\nc}\overline\R_{g,\ell}^0$
coincides with the $T$-locus for these values of $\ell$.
\end{rmk}

To show Theorem \ref{teo_jt} we will prove a stronger proposition.
\begin{prop}
Given a rooted curve $(\ssC,\ssL,\phi)$ of genus $g\geq 4$ which
is not a $J$-curve, if $\msf a\in\underline\Aut(\ssC,\ssL,\phi)\slash\QR(\ssC,\ssL,\phi)$
is a junior automorphism, then its coarsening $a$ is an elliptic tail automorphism 
of order $3$ or $6$.
\end{prop}
\proof
We introduce the notion of $\star$-smoothing, following
\cite{harmum82} and~\cite{lud10}. 
\begin{defin}
Consider a rooted curve $(\ssC,\ssL,\phi)$,
suppose $\msf{a}\in\underline\Aut(\ssC,\ssL,\phi)$ is an automorphisms
such that there exists a cycle of $m$ non-separating nodes $\msf{n}_0,\dots,\msf{n}_{m-1}$,
\cor{i.e.}~we have $\msf a(\msf n_i)=\msf n_{i+1}$ for all $i=0,1\dots,m-2$ and $\msf{a} (\msf{n}_{m-1})=\msf n_0$.
The triple is $\star$-smoothable if and only if the action of $\msf a^m$ over
the coordinate associated to every node is trivial. This is equivalent to ask $\msf a^m(t_{\msf n_i})=t_{\msf n_i}$
for any $i=0,1,\dots,m-2$, where $t_{\msf n_i}$ is the coordinate associated
to node~$\msf n_i$.
\end{defin}

If $(\ssC,\ssL,\phi)$ is $\star$-smoothable, there exists a deformation of $(\ssC,\ssL,\phi)$
together with an automorphism $\msf a$ and 
smoothing the $m$ nodes of the cycle.
Moreover, this deformation preserves the age of the $\msf a $-action on $\defo(\ssC,\ssL,\phi)\slash\QR$. Indeed,
the eigenvalues of $\msf a$ are a discrete and locally constant set, thus constant by
deformation.
As the $T$-locus and the $J$-locus are closed by $\star$-smoothing, we can suppose,
as an additional hypothesis of Theorem \ref{teo_jt}, that our curves are $\star$-\cor{rigid}, \cor{i.e.}~non-$\star$-smoothable.
From this point we suppose that $(\ssC,\ssL,\phi)$ is $\star$-rigid.

We will show in eight steps that if the group
$$\underline\Aut(\ssC,\ssL,\phi)\slash\QR\left(\ssC,\ssL,\phi\right)$$
is junior, and $(\ssC,\ssL,\phi)$ is not a $J$-curve, then it is a $T$-curve.
After the age criterion \ref{age_crit} and Proposition \ref{propt}, this will prove
Theorem \ref{teo_jt}. 
From now on we work under the hypothesis that $\msf a\in\underline\Aut_C(\ssC,\ssL,\phi)\slash\QR$
is a non-trivial automorphism aged less than $1$, and $(\ssC,\ssL,\phi)$ is not a $J$-curve.\newline

\cor{Step 1.}
Consider the decorated graph $(\Gamma_0,M)$ of $(\ssC,\ssL,\phi)$.
As before, we call $E_{\sep}$ the set of separating edges of $\Gamma_0$.
Then, following Remark \ref{defo}, we can estimate the age by a splitting of the form
$$\bigoplus_{e\in E_{\sep}}\A_{t_e}\oplus \bigoplus_{e'\in E\backslash E_{\sep}}\A_{t_{e'}}
\oplus\bigoplus_{i=1}^mH^1(\overline\ssC_i,T_{\overline\ssC_i}(-D_i)),$$
where $t_e$ is a coordinate parametrizing the smoothing of the node associated to edge $e$,
and the curves $\overline\ssC_i$ are the normalizations of the irreducible components of $\ssC$.

Every automorphism in $\underline\Aut(\ssC,\ssL,\phi)$ fixes the three summand
in the sum above. Moreover, every quasireflection acts only on the first summand,
then, by Propostion \ref{propt}, the group $\underline\Aut(\ssC,\ssL,\phi)\slash\QR$ acts on
\begin{equation}\label{eq_space}
\left( \frac{\bigoplus_{e\in E_{\sep}}\A_{t_e}}{\QR(\ssC,\ssL,\phi)}\right)\oplus
 \bigoplus_{e'\in E\backslash E_{\sep}}\A_{t_{e'}}
\oplus\bigoplus_{i=1}^mH^1(\overline\ssC_i,T_{\overline\ssC_i}(-D_i)).
\end{equation}

Every quasireflection acts on exactly one coordinate $t_e$
with $e\in E_{\sep}$. We rescale all the coordinates $t_e$
by the action of $\QR(\ssC,\ssL,\phi)$. We call $\tilde t_e$,
for $e\in E(\Gamma_0)$, the new set of coordinates.
Obviously $\tilde t_{e'}=t_{e'}$ if $e'\in E(\Gamma_0)\backslash E_{\sep}$.\newline

\cor{Step 2.} We show two lemmata about the age contribution of
the $\msf a$-action on nodes, which we will call \emph{aging} on nodes.

\begin{defin}[coarsening order]
If $\msf a\in\underline\Aut(\ssC,\ssL,\phi)$ and $a$ is its coarsening, then we define 
$$\cord\msf a:=\ord a.$$
\end{defin}
The coarsening order is the least integer $m$ for which $\msf a^m$ is a ghost automorphism.

\begin{lemma}\label{lemma_bound}
Suppose that $\msf Z\subset \ssC$ is a subcurve of $\ssC$ such that $\msf a(\msf Z)=\msf Z$ and 
$\msf n_0,\dots \msf n_{m-1}$ is a cycle, by $\msf a $, of nodes in $\msf Z$. Then
we have the following inequalities:
\begin{enumerate}
\item $\age(\msf a)\geq \frac{m-1}{2}$;
\item  if the nodes are non-separating, $\age (\msf a)\geq \frac{m}{\ord\left(\msf a|_{\msf Z}\right)}+\frac{m-1}{2}$;
\item if $\msf a^{\cord \msf a}$ is a senior ghost, we have $\age (\msf a)\geq \frac{1}{\cord(\msf{a})}+\frac{m-1}{2}$.
\end{enumerate}
\end{lemma}
\proof
We call $\tilde t_0,\tilde t_1,\dots,\tilde t_{m-1}$ the coordinates associated respectively to nodes
$\msf n_0,\dots,\msf n_{m-1}$. By hypothesis, $\msf a (\tilde t_0)=c_1\cdot \tilde t_1$ and
$\msf a^i(\tilde t_0)=c_i\cdot\tilde t_i$ for all $i=2,\dots,m-1$,
where $c_i$ are complex numbers. If $n=\ord(\msf a|_{\msf Z})$,
we have $\msf a^m(\tilde t_0)=\xi_n^{u m}\cdot \tilde t_0$ where $\xi_n$
is the primitive $n$th root of unity and $0 \leq u <n\slash m$.
We
call $u$ exponent of the cycle.
Observe that $\msf a(\tilde t_{i-1})=(c_i\slash c_{i-1})\cdot \tilde t_i$ and $\msf a^m (\tilde t_i)=\xi_n^{u m}\cdot\tilde t_i$
for all $i$.

We can explicitly write the eigenvectors for the action of $\msf a$ on the coordinates
$\tilde t_0,\dots,\tilde t_{m-1}$. Set $d:=n\slash m$ and $b:=sd+u$ with $0\leq s<m$, and consider the vector
$$v_b:=(\tilde t_0=1,\ \tilde t_1=c_1\cdot \xi_n^{-b},\dots,\ \tilde t_i=c_i\cdot\xi_n^{-ib},\dots).$$
Then $\msf a (v_b)=\xi_n^b\cdot v_b$. The contribution to the age of the eigenvalue $\xi_n^b$ is
$b\slash n$, thus we have
$$\age\msf a = \sum_{s=0}^{m-1}\frac{sd+u}{n}=\frac{mu}{n}+\frac{m-1}{2},$$
proving point $(1)$.

If the nodes are non-separating, as we are working on a $\star$-rigid curve,
we have $u\geq 1$ and the point $(2)$ is proved.

Suppose that the action of $\msf a$ on $\ssC$ has $j$ nodes cycles
of order $m_1, m_2,\dots, m_j$, of exponents respectively $u_1,\dots, u_j$ . If $k=\cord \msf a$,
 $\msf a^k$ fixes every node and
its age is $(\sum m_iu_ik)\slash n$, which is greater
or equal to $1$ by hypothesis.
By the previous result, the age of $\msf a$ on the $i$th cycle
is bounded by $m_iu_i\slash n+(m_i-1)\slash 2$.
As a consequence
$$\age \msf a= \sum_{i=1}^j\left(\frac{m_iu_i}{n}+\frac{m_i-1}{2}\right)\geq \frac{1}{k}+\frac{m_1-1}{2}.$$ \fine

\begin{lemma}\label{lemma_node}
Suppose that $(\ssC,\ssL,\phi)$ is $\star$-rigid and 
consider a non-separating node $\msf n$. If
the coarsening of
$\msf a$ 
acts locally exchanging the branches of $\msf n$,
then the $\msf a$-action yields an aging of $1\slash2 $
on the nodes.
\end{lemma}
\proof
The automorphism $\msf a$ induces an obvious
automorphism of the decorated graph $(\Gamma_0,M)$.
We will call this automorphism $\msf a$ too, with a little abuse of notation.
Thus we have $\msf a^*M=M$. In particular, if $e_{\msf n}$
is an oriented edge associated to $\msf n$, then 
$M(e_{\msf n})\equiv \msf a^*M(e_{\msf n})\equiv M(\overline e_{\msf n})\equiv -M(e_{\msf n})$.
Therefore $M(e_{\msf n})\equiv\ell\slash 2$.
If $\ell$ is prime thus necessarily $\ell=2$. Anyway, also
when $\ell$ is composite, by the generalized definition of $M$ (see Remark \ref{rmk_composite})
we have that the order of the local stabilizer is $r=2$ and $\msf a(e_{\msf n})\in\Z\slash 2$.
As $(\ssC,\ssL,\phi)$ is $\star$-rigid, $\msf a (e_{\msf n})\not\equiv 0$,
therefore the aging is~$1\slash 2$.\fine

\cor{Step 3.} We observe that all the nodes of $\ssC$ are fixed except at most two of
them, who are exchanged. Moreover, if a pair of non-fixed
nodes exists, they contribute by at least $1\slash 2$.
This fact is a straightforward consequence of the first
point in Lemma \ref{lemma_bound}.\newline

\cor{Step 4.} Consider an irreducible component $\msf Z\subset \ssC$, then $\msf a(\msf Z)=\msf Z$.
To show this, we suppose there exists a cycle of irreducible components $\msf C_1, \dots,\msf C_m$  with $m\geq 2$ such that
$\msf a(\msf C_i)=\msf C_{i+1}$ for $i=1,\dots,m-1$, and $\msf a(\msf C_m)=\msf C_1$.
We call $\overline\ssC_i$ the normalizations of these components, and
$D_i$ the preimages of $\ssC$ nodes on $\overline\ssC_i$.
We point out that this construction implies that $(\overline\ssC_i,D_i)\cong (\overline\ssC_j, D_j)$ for all $i,j$.
Then,
an argument of \cite[p.34]{harmum82} shows that the action of $\msf a$ on
$$ \bigoplus_{i=1}^mH^1(\overline\ssC_i, T_{\overline\ssC_i}(-D_i))\subset \defo(\ssC)$$
gives a contribution of at least $k\cdot(m-1)\slash2$ to $\age\msf a$, where
$$k=\dim H^1(\overline\ssC_i,T_{\overline\ssC_i}(-D_i))=3g_i-3+\#D_i.$$

This give us two cases for which $m$ could be greater than $1$
with still a junior age: $k=1$ and $m=2$ or $k=0$.

If $k=1$ and $m=2$, we have $g_i=0$ or $1$ for $i=1,2$.
Moreover, the aging of at least $1\slash 2$ sums to
another aging of $1\slash 2$ if there is a pair of non-fixed nodes.
As $\msf a$ is junior, we conclude that $\ssC=\ssC_1\cup\msf a(\ssC_1)$
but this implies $g(\ssC)\leq 3$, contradiction.

If $k=0$, we have $g_i=1$ or $g_i=0$, the first is excluded because
it implies $\#D_i=0$ but the component must intersect the curve somewhere.
Thus, for every component in the cycle, the normalization $\overline\ssC_i$ is the
projective line $\p^1$ with $3$ marked points. We have two cases:
the component $\ssC_i$ intersect $\overline{\ssC\backslash\ssC_i}$ in $3$ points
or in $1$ point, in the second case $\ssC_i$ has an autointersection node
and $\ssC=\ssC_1\cup\msf a(\ssC_1)$, which is a contradiction because $g(\ssC)\geq 4$.
It remains the case in the image below.

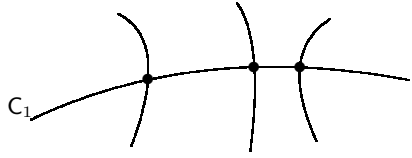
\begin{figure}[h]
 \begin{picture}(400,50)(40,0)
   \qbezier(187,10)(250,40)(330,25)
   \qbezier(220,50)(240,40)(225,0)
   \qbezier(265,54)(275,40)(270,-2)
   \qbezier(300,48)(280,35)(295,2)
   
   \put(178,12){\m{\footnotesize{$ \ssC_1$}}}
   \put(231.2,25.3){\circle*{4}}
   \put(271.3,30){\circle*{4}}
   \put(288.8,30){\circle*{4}}
 
 \end{picture}
 \caption{Case with $ \ssC_1\cong \p^1$ and $3$ marked points}
\end{figure} 

As $\ssC_1,\ssC_2,\dots,\ssC_m$ are moved by $\msf a$, every
node on $\ssC_1$ is transposed with another one or is fixed 
with its branches interchanged. In both cases the aging is $1\slash 2$,
after a straightforward analysis we obtain an age contribution bigger than $1$ using Lemmata \ref{lemma_bound}
and \ref{lemma_node}.\newline

\cor{Step 5.} We prove that every node is fixed by $\msf a$.
Consider the normalization $\msf{nor}\colon\bigsqcup_i\overline\ssC_i\to\ssC$
already introduced.
If the age of $\msf a$ is lower than $1$, \emph{a fortiori} 
we have $\age \msf a|_{\overline\ssC_i}<1$ for all $i$.
In \cite[p.28]{harmum82}
there is a list of those smooth stable curves for which there exists
a non-trivial junior action.

\begin{enumerate}[label=\roman{*}., ref=(\roman{*})]
\item The projective line $\p^1$ with $\msf a\colon z\mapsto (-z)$ or $(\xi_4 z)$;
\item an elliptic curve with $\msf a$ of order $2,\ 3,\ 4$ or $6$;
\item an hyperelliptic curve of genus $2$ or $3$ with $\msf a$ the hyperelliptic involution;
\item a bielliptic curve of genus $2$ with $\msf a$ the canonical involution.
\end{enumerate}
We observe that the order of the $\msf a$-action on these components is always $2,\ 3,\ 4$ or $6$.
As a consequence, if $\msf a$ is junior, then $n=\cord\msf a=2,3,4,6$ or~$12$,
as it is the greatest common divisor between the $\cord\left(\msf a|_{\overline\ssC_i}\right)$.

First we suppose $\ord\msf a>\cord\msf a$, thus $\msf a^{\cord\msf a}$ is a ghost and it must
be senior. Indeed, if $\msf a^{\cord\msf a}$ is aged less than $1$,
then $(\ssC,\ssL,\phi)$ admits junior ghosts,
contradicting our assumption. By point $(3)$ of
Lemma \ref{lemma_bound}, if there
exists a pair of non-fixed nodes, we obtain an aging of $1\slash n+1\slash2$ on
node coordinates.
If  $\ord \msf a=\cord\msf a$ the bound is even greater.
As every component is fixed by $\msf a$, the two nodes are non-separating, and by point $(2)$
of Lemma \ref{lemma_bound} we obtain an aging of $2\slash n+1\slash2$.

If $\overline\ssC_i$ admits an automorphism of order $3,\ 4$ or $6$, 
by a previous 
analysis of Harris and Mumford (see \cite{harmum82} again),
this yields
an aging  of, respectively, $1\slash 3$,\ $1\slash 2$ and $1\slash 3$ on $H^1(\overline\ssC_i,T_{\overline\ssC_i}(-D_i))$.

These results combined, show that a non-fixed pair of nodes gives
an age greater than $1$. Thus, if $\msf a$ is junior,
every node is fixed.\newline

\cor{Step 6.} We study the action of $\msf a$ separately on every irreducible
component. The $\msf a$-action is non-trivial on at least one component~$\ssC_i$, and this component must lie
in the list above. 

In case (i), $\overline\ssC_i$ has at least $3$ marked points
because of the stability condition. Actions of type $x\mapsto \zeta x$
have two fixed points on $\p^1$, thus
at least one of the marked points is non-fixed.
A non-fixed preimage of a node has order $2$, thus
 the coarsening $a$ of $\msf a$ is the involution $z\mapsto -z$. Moreover,
$\ssC_i$ is the autointersection of the projective line and
$\msf a$ exchanges the branches of the node. Because of Lemma \ref{lemma_node},
$\msf a$ acts non-trivially with order $2$ on the node, and its
action on the associated coordinate gives an aging of $1\slash2$.

The analysis for cases (iii) and (iv) is identical to that
developed in \cite{harmum82}: the only possibilities of a junior action is the case
of an hyperelliptic curve $\msf E$ of genus $2$ intersecting $\overline{\ssC\backslash\msf E} $
in exactly one point, whose hyperinvolution gives an 
aging of $1\slash 2$ on $H^1(\overline\ssC_i,T_{\overline\ssC_i}(-D_i))$.

Finally, in case (ii), we use again the analysis of \cite{harmum82}.
The elliptic component $\msf E$ has
$1$ or $2$ point of intersection with $\overline{\ssC\backslash\msf E}$.
If there is $1$ point of intersection, elliptic tail case,
for a good choice of coordinates
the coarsening $a$ acts as $z\mapsto \xi_nz$,
where $n$ is $2,3,4$ or $6$. The aging is, respectively,
$0,\ 1\slash3,\ 1\slash2,\ 1\slash3.$ If there are $2$ points of
intersection, elliptic ladder case, the order of $\msf a$ on $\msf E$ must be $2$ or $4$ and the
aging respectively $1\slash 2$ or $3\slash 4$.\newline

\cor{Step 7.} 
Resuming what we saw until now, if $\msf a$ is a junior
automorphism of $(\ssC,\ssL,\phi)$, $a$ its coarsening
and $C_1$ an irreducible component of $C$,
then we have one of the following:
\begin{enumerate}
\item [A.] component $C_1$ is an hyperelliptic tail, crossing the curve in one point,
with $a$ acting as the hyperelliptic involution and aging $1\slash 2$ on $H^1(\overline C_1,T_{\overline\ssC_i}(-D_1))$;
\item [B.] component $C_1$ is a projective line $\p^1$ autointersecting itself, crossing the curve in one point,
with $a$ the involution which fixes the nodes, and aging $1\slash 2$;
\item [C.] component $C_1$ is an elliptic ladder, crossing the curve in two points,
with $a$ of order $2$ or $4$ and aging respectively $1\slash2$ or $3\slash 4$;
\item [D.] component $C_1$ is an elliptic tail, crossing the curve
in one point, with $a$ of order $2,3,4$ or $6$ and aging $0,1\slash3,1\slash 2$ or $1\slash 3$;
\item [E.] automorphism $a$ acts trivially on $C_1$ with no aging.
\end{enumerate}

\begin{figure}[h]
 \begin{picture}(400,70)(40,0)
 \put(20,-10){
 \qbezier(20,70)(50,60)(90,90)
 \qbezier(70,90)(80,65)(75,10)
 
 \put(73.4,79.3){\circle*{3}}
 
 \put(-1,76){\m{\footnotesize{$g(C_1)$$=$$2$}}}
 
 \put(67.5,69.4){\m{\footnotesize{$\msf n$}}}}
 
 \put(-20,0){
 \qbezier(187,30)(210,50)(280,40)
 \qbezier(280,40)(315,35)(295,15)
 \qbezier(295,15)(290,10)(285,15)
 \qbezier(285,15)(275,25)(290,60)
 
 \qbezier(210,70)(230,40)(215,10)
 
 \put(221,42.4){\circle*{3}}
 \put(211,46){\m{\footnotesize{$\msf n$}}}
 
 \put(282.8,39.9){\circle*{3}}

\put(240,48){\m{\footnotesize{$\overline C_1=\p^1$}}} }

 \put(335,-10){
 \qbezier(0,80)(50,60)(120,90)
 \qbezier(20,90)(30,65)(25,10)
 \qbezier(90,90)(80,65)(95,10)
 
 \put(87,78.6){\circle*{3}}
 \put(24.4,74){\circle*{3}}
 
 \put(31,78){\m{\footnotesize{$g(C_1)$$=$$1$}}}}
 
 \end{picture}
 \caption{Components of type A, B and C.}
\end{figure}
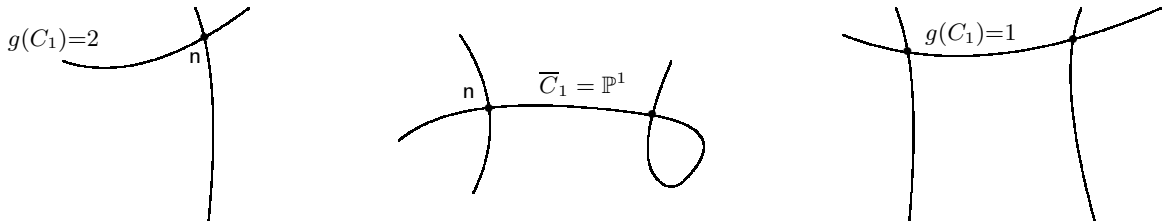 

We rule out cases (A), (B) and (C).
At first we suppose
there is a component of type (A) or (B). For genus reasons, the component
intersected in both cases must be of type (E).
We study the local action on
the separating node $\msf n$. The local picture of $\msf n$ is 
$\left[\{xy=0\}\slash\mmu_r\right]$
where $r$ is the order of the local stabilizer.
The smoothing
of the node is given by the stack $\left[\{xy=t_{\msf{n}}\}\slash\mmu_r\right]$.
We observe that there
exists a quasireflection of order $r$ acting non-trivially on the node,
and it generates quasireflections acting on this node.
Thus $\tilde t_{\msf n}=t_{\msf n}^r$.
For a good choice of coordinates we have that the local picture of the
coarse space $C$ at the node is $\{x'y'=0\}$, with $x'=x^r$
and $y'=y^r$. Therefore the local deformation of the coarse
space is given by $\{x'y'=\tilde t_{\msf n}\}$. The local action of $\msf a$ on the coarse space
is $x'\mapsto -x'$ and $y'\mapsto y'$, thus the action of
$\msf a$ lifts to $\tilde t_{\msf n}\mapsto -\tilde t_{\msf n}$. The
additional age contribution is $1\slash 2$, ruling out this case.

In case there is a component of type (C), if its nodes
are separating, then one of them must intersect a component
of type (E) and we use the previous idea. In case nodes are non-separating, 
we use Lemma \ref{lemma_bound}. If $\ord \msf a>\cord \msf a$, then $\msf a^{\cord\msf a}$ is
a senior ghost because $(\ssC,\ssL,\phi)$ is not a $J$-curve, thus by
point $(3)$ of the lemma there is an aging of $1\slash\cord\msf a$ on the node coordinates.
If $\ord\msf a=\cord\msf a$, the bound is even greater, as by point $(2)$
we have an aging of $2\slash\cord\msf a$. We observe that $\cord\msf a=2,4$ or $6$,
and in case $\cord\msf a=6$ there must be a component of type (E).
Using additional contributions listed above we rule out
the case (C).\newline

\cor{Step 8.} We proved that $\ssC$ contains
components of type (D) or (E), \cor{i.e.}~the automorphism
$\msf a$ acts non-trivially only on elliptic tails.
If $\msf n$ is the elliptic tail node, there are
possibly two quasireflections acting on the coordinate $t_{\msf n}$:
a ghost automorphism associated to this node
and the elliptic tail quasireflection. If the order of the local stabilizer is $r$,
then $\tilde t_{\msf n}=t_{\msf n}^{2r}$.

If $\ord \msf a=2$ we are in the ETQR case, this action is a quasireflection
and it contributes to rescaling the coordinate $t_{\msf n}$.

If $\ord\msf a=4$, the action on the (coarse) elliptic tail is $z\mapsto \xi_4z$. The space
$H^1(\overline\ssC_i,T_{\overline\ssC_i}(-D_i))$ is the space of $2$-forms
$H^0(\overline\ssC_i,2K_i)$: this space is generated by $dz^{\xx 2}$
and the action of $\msf a$ is $dz^{\xx 2}\mapsto\xi_4^2dz^{\xx2}$.
Moreover, if the local picture of the elliptic tail node is $[\{xy=0\}\slash\mmu_r]$,
then $\msf a\colon(x,y)\mapsto(\zeta x,\zeta' y)$ such that
$\zeta^r=\xi_4$ and $(\zeta')^r=1$. As a consequence $\msf a\colon t_{\msf n}\mapsto \zeta\zeta' t_{\msf n}$
and therefore $\tilde t_{\msf n}\mapsto \xi_2\tilde t_{\msf n}$.
Then, $\age \msf a=1\slash2+1\slash2$, proving the seniority of $\msf a$.

If $\msf E$ admits an automorphism $\msf a$ of order $3$ or $6$,
 the action
on the (coarse) elliptic tail is $\msf a\colon z\mapsto \xi_6^kz$.
Then $dz^{\xx 2}\mapsto \xi_3^kdz^{\xx2} $ and $\tilde t_{\msf n}\mapsto\xi_3^k\tilde t_{\msf n}$.
For $k=1,4$ we have age lower than $1$.\newline

If $(\ssC,\ssL,\phi)$ is not a $J$-curve, 
we have shown that
 the only case where $\msf a\in\underline\Aut(\ssC,\ssL,\phi)\slash\QR$ is junior,
is when its coarsening $a$ is an elliptic tail automorphism of order $3$ or $6$.\fine

The previous theorem reduces the analysis of $\sing^{\nc}\overline\R_{g,\ell}^k $
to two loci. The next section will be devoted to the $J$-locus.
About the $T$-locus we observe that it is
that part of the non-canonical locus ``coming'' from $\overline\M_g$.
More formally, if $\pi\colon\overline\R_{g,\ell}^k\to\overline\M_g$ is
the natural projection, we have
$T_{g,\ell}^k\subset\pi^{-1}\sing^{\nc}\overline\M_g=T_{g,1,0}.$
Harris and Mumford showed that 
multicanonical forms extend over the $T$-locus
on $\overline\M_g$. Their
proof could be adapted for $\overline\R_{g,\ell}^k$ 
as shown precisely in \cite[Theorem 4.1]{lud10}
for the case $\ell=2$, $k=1$.

In the case of the $J$-locus it is unknown if the extension
of the pluricanonical forms is possible. Harris and Mumford's techniques
do not adapt and there is at the moment no global approach 
to treat this case. In the next section we will give a description of the $J$-locus
based on the dual graph structure. In particular we will show the cases
with non-empty $J$-locus for small values of $\ell$.

\section{The $J$-locus}

The $J$-locus is the ``new'' part of the non-canonical locus
which appears passing from $\overline\mm_g$ to
one of its coverings $\overline\rr_{g,\ell}^k$.  
For some values of $\ell$ and $k$ this locus could
 be empty, however, a consequence of our analysis
is that it is actually not empty for any $\ell>2$
and $k\neq 0$.

We will exhibit an explicit decomposition of $J_{g,\ell}^k$
in terms of the strata $\s_{(\Gamma,M)}$.
We point out a significant difference with respect to the
description of the singular locus: in the case of $H_{g,\ell}^k$,
that we obtained in Section \ref{subs_strata},
we showed a decomposition in terms of loci whose generic
point represents a two component curve, \cor{i.e.}~a
vine curve. For the $J$-locus we do not have
such an elegant minimal decomposition: for values
of $\ell$ large enough, there exists strata representing 
$J$-curves with an arbitrary high number of components,
such that each one of their smoothing is not a $J$-curve.
Equivalently,
there are decorated graphs with an arbitrary high
number of vertices and admitting junior
automorphisms, but such that each one of its contraction
does not admit junior automorphisms.

\begin{rmk}
We already observed that a ghost automorphism
always acts trivially on loop edges of decorated dual graphs,
and that quasireflections only act on separating edges.
Thus we can ignore these edges in studying
$\underline\Aut_C(\ssC,\ssL,\phi)\slash\QR_C(\ssC,\ssL,\phi)$.
From this point we will automatically contract
loops and separating edges as they appear.
This is not a big change in our setting,
in fact graphs without loops and separating edges are a subset of the graphs we considered
until now. Reducing our analysis to this subset is done for the
sole purpose of simplifying the notation.
\end{rmk}

Age is not well-behaved with respect to
graph contraction, but there is another invariant which is
better behaved: we will define a number associated
to every ghost, which
respects a super-additive property in the case
of strata intersection (see Theorem \ref{teo_levage}). Here we only work under the condition that $\ell$ be a prime number.

At first consider a decorated graph contraction
$$(\Gamma_0,M)\to (\Gamma_1,M_1).$$
From the definition of the stratification, we know that this implies
$$\s_{(\Gamma_1,M_1)}\subset \overline\s_{(\Gamma_0,M)}.$$
We know that $\underline\Aut_C(\ssC,\ssL,\phi)$
is canonically isomorphic to the group 
$G(\Gamma_0;\Z\slash\ell)\subset S(\Gamma_0;\Z\slash\ell)$
for every rooted curve $(\ssC,\ssL,\phi)$ in $\s_{(\Gamma_0,M)}$,
and moreover $S(\Gamma_0;\Z\slash\ell)$ is canonically identified with $C^1(\Gamma_0;\Z\slash\ell)$
via $M$ multiplication.
 As there is no risk of confusion, from now on
we will not repeat the notation of $\Z\slash\ell$ on the label of groups $G$, $S$ and $C^1$.

Because of contraction, the edge set $E(\Gamma_1)$ is a subset of $E(\Gamma_0)$,
so that the group $S(\Gamma_1)$ is naturally immersed in $S(\Gamma_0)$,
and $C^1(\Gamma_1)$ in $C^1(\Gamma_0)$.
These immersions are compatible with multiplication by $M$, 
thus $G(\Gamma_1)=G(\Gamma_0)\cap S(\Gamma_1)$ and
we have a natural immersion of the $G$ groups too.
\begin{prop}\label{prop_H1}
Given a contraction $(\Gamma_0,M)\to(\Gamma_1,M_1)$,
the elements of $G(\Gamma_1)$ are those cochains of $G(\Gamma_0)$ whose
support is contained in $E(\Gamma_1)$.
\end{prop}
This gives an interesting correspondence between curve specialization,
decorated graph contraction, strata inclusion and canonical immersions between
the associated ghost automorphism groups.

From the definition of $G(\Gamma_0;\Z\slash\ell)$ and Proposition \ref{cut_gen},
we know that $G(\Gamma_0)\cong\Z\slash\ell^{\#V(\Gamma_0)-1}$,
and we have an explicit basis for it. We consider a spanning tree T for $\Gamma_0$,
we call $e_1, e_2,\dots, e_k$ the edges of $T$, each one with an orientation, such
that $k=\#V(\Gamma_0)-1$. Then the cuts $\cut_{\Gamma_0}(e_i;T)$ form a basis
of $G(\Gamma_0)$. We can also write
$$G(\Gamma_0)=\bigoplus_{i=1}^k\left(\cut_{\Gamma_0}(e_i;T)\cdot\Z\slash\ell\right).$$
For an $n$-vine decorated graph $(\Gamma_{(2,n)},M)$, the
$G$-group is  cyclic.

\begin{rmk}
An element $\msf a\in G(\Gamma_0)$ could be seen
as living on stratum $\s_{(\Gamma_0,M)}$. We observe that
if $(\Gamma_0,M)\to(\Gamma_1,M_1)$ is a contraction,
then $\msf a$ lies in $G(\Gamma_1)$ by
the natural injection. Thus every automorphism living on $\s_{(\Gamma_0,M)}$
lives on all the closure $\overline\s_{(\Gamma_0,M)}$.
\end{rmk}

If a ghost automorphism is junior, it carries a non-canonical
singularity that spreads all over the closure of the stratum where the automorphism lives.
This informal statement justifies the following definition.

\begin{defin}
The age of a stratum $\s_{(\Gamma_0,M)}$ is the minimum age of a ghost automorphism
$\msf a$ in $G(\Gamma_0)$. 
As in the case of group age, the age of an automorphism depends on the
primitive root chosen, but the stratum age does not.
\end{defin}

With this new notation, the locus of non-canonical singularities could be written as follows,
$$\sing^{\nc}\overline\R_{g,\ell}^k=\bigcup_{\age\s_{(\Gamma_0,M)}<1}\overline\s_{(\Gamma_0,M)}.$$
Indeed, if $[\ssC,\ssL,\phi]\in\overline\R_{g,\ell}^k$ has a junior ghosts group, then
this point lies on the closure of a junior stratum. Conversely, every point in the closure of
a junior stratum has a junior ghosts group.

\begin{defin}
We say that a set of contractions $\{(\Gamma_0,M)\to(\Gamma_i,M_i)\}$, for
$i=1,\dots,k$, \cor{covers} $\Gamma_0$ if $E(\Gamma_0)=\bigcup_iE(\Gamma_i)$.
We observe that if this set is a covering, then $G(\Gamma_i,M_i)\subset G(\Gamma_0,M)$
for all $i$ and
$$\overline\s_{(\Gamma_0,M)}\subset\bigcap_{i=1}^k\overline\s_{(\Gamma_i,M_i)}.$$
\end{defin}

\begin{defin}
An automorphism $\msf a\in G(\Gamma_0)$ is \cor{supported} on
stratum $\s_{(\Gamma_0,M)}$ if its support is the whole $E(\Gamma_0)$ set.
We observe that this property has an immediate moduli interpretation:
an automorphism supported on $\s_{(\Gamma_0,M)}$ appears in the ghost group
of every curve in $\overline\s_{(\Gamma_0,M)}$ but it does not appear
in any other stratum whose closure contains $\s_{(\Gamma_0,M)}$. 
\end{defin}

Unfortunately, the age of a strata intersection is not bounded by the sum of the ages
of strata intersecting. However, there exists another invariant which has a superadditive
property with respect to strata intersection. We will pay attention to the new
automorphisms that appear at the intersection, \cor{i.e.}~those automorphisms
supported on the intersection stratum, using the notion just introduced.

\begin{teo}\label{teo_levage}
Consider a covering $(\Gamma_0,M)\to(\Gamma_i,M_i)$, with $i=1,\dots, m$,
such that 
$$G(\Gamma_0)=\sum_{i=1}^mG(\Gamma_i).$$
Then for every $\msf a$ supported on $\s_{(\Gamma_0,M)}$
we have
$$\age \msf a-\#E(\Gamma_0)\geq \sum_{i=1}^m\left(\age\s_{(\Gamma_i,M_i)}-\#E(\Gamma_i)\right).$$
\end{teo}

To prove the theorem we need the following lemma

\begin{lemma}\label{lem_tec2}
If $\msf{a}$ is supported on $\s_{(\Gamma_0,M)}$, then 
$$\age\msf{a}+\age\msf{a}^{-1}=\#E(\Gamma_0).$$
\end{lemma}
\proof
Given an edge $e\in E(\Gamma_0)$,
by definition, $\msf{a}^{-1}(e),\msf{a}(e)\in\Z\slash\ell$. As $\msf{a}$ is supported on $\s_{(\Gamma,M)}$,
$\msf{a}(e)\not\equiv 0$ for all $e$ in $E(\Gamma_0)$, and then this component brings
to $\age\msf{a}$ and $\age\msf{a}^{-1}$ respectively a value of $\msf{a}(e)\slash\ell$ and $\left(1-\msf{a}(e)\slash\ell\right)$.
As a consequence we obtain $\age\msf{a}+\age\msf{a}^{-1}=\#E(\Gamma)=\codim\s_{(\Gamma_0,M)}$.\fine

As a direct consequence of the previous lemma,
we have $\age\msf a^{-1}=\#E(\Gamma_0)-\age\msf a$.
By hypothesis we can write,
$$\msf a^{-1}=\msf a_1+\msf a_2+\cdots+\msf a_m,$$
where $\msf a_i\in G(\Gamma_i,M_i)$ for all $i$, and we call
$c_i$ the cardinality of $\msf a_i$ support. By subadditivity of age, we have
$\age \msf a^{-1}\leq \sum\age\msf a_i$, then using Lemma \ref{lem_tec2}
 we obtain
$$\msf a-\#E\geq \sum_{i=1}^m\left(\age \msf a_i^{-1}-c_i\right).$$
By the fact that $c_i\leq \#E(\Gamma_i)$ for all $i$,
and by the definition of age for the strata, the Theorem is proved.\fine

We observe that $\#E(\Gamma_i)=\codim\s_{(\Gamma_i,M_i)}$ by 
Proposition \ref{prop_codim}, thus we found an inequality about
 age involving geometric data. There is another formulation
of the statement. We already observed that for every graph $\Gamma$,
$b_1(\Gamma)=\#E-(\#V-1)$. If the sum in the hypothesis is a direct sum,
the rank condition is equivalent to $\#V(\Gamma_0)-1=\sum \#V(\Gamma_i)-1$.
Therefore we have the following.
\begin{coro}
If the $(\Gamma_i,M_i)$ cover $(\Gamma_0,M)$ and
$$G(\Gamma_0)=\bigoplus_{i=1}^m G(\Gamma_i),$$
then for every automorphism $\msf a$
supported on $G(\Gamma_0)$ we have
$$\age\msf a-b_1(\Gamma_0)\geq\sum_{i=1}^m\left(\age\s_{(\Gamma_i,M_i)}-b_1(\Gamma_i)\right).$$
\end{coro}

In the case of two strata intersecting, a rank condition
implies the splitting of $G(\Gamma_0,M)$ in a direct sum.

\begin{lemma}\label{lem_cover}
Consider the contractions of decorated graphs $(\Gamma_0,M)\to (\Gamma_i,M_i)$
 for $i=1,2$. If $\Gamma_1$ and $\Gamma_2$ cover $\Gamma_0$, and moreover
$$\#V(\Gamma_0)-1=(\#V(\Gamma_1)-1)+(\#V(\Gamma_2)-1),$$
then we have
$$G(\Gamma_0)= G(\Gamma_1)\oplus G(\Gamma_2).$$
\end{lemma}
\proof
Before proving it, we point out a useful fact:
given any contraction
$\Gamma_0\to\Gamma_i$, the natural injection
$G(\Gamma_i)\hookrightarrow G(\Gamma_0)$
sends cuts on cuts.
We suppose, without loss of generalities, that 
$\#V(\Gamma_1)\leq \#V(\Gamma_2)$
and we prove the lemma by induction on $\#V(\Gamma_1)$.

The base case $\#V(\Gamma_1)=1$ is empty,
$\Gamma_0=\Gamma_2$ and the thesis follows obviously.
Now suppose $\#V(\Gamma_1)=q>1$, then $\#V(\Gamma_2)<\#V(\Gamma_0)$
and so there exists two vertices of $\Gamma_0$ connected by edges who lies
in $E(\Gamma_1)$ but not in $E(\Gamma_2)$. We call $e_1$ one of these
edges in $E(\Gamma_1)$ and $T_1$ a spanning tree of $\Gamma_1$ containing $e_1$. If 
$\cut_{\Gamma_1}(e_1;T_1)$ is the corresponding cut, it is also
an element of $G(\Gamma_0)$.
By rank conditions it suffices to prove that $G(\Gamma_0)= G(\Gamma_1)+G(\Gamma_2)$.
Consider $\msf a\in G(\Gamma_0)$ which is not a sum of elements
in $G(\Gamma_1)$ and $G(\Gamma_2)$. Now we define
$$\msf a':=\msf a - k\cdot \cut_{\Gamma_1}(e_1;T_1),$$
where $k$ is the necessary integer such that $\msf a'(e_1)\equiv0$.
Consider the graphs $(\Gamma_0',M')$ and $(\Gamma_1',M_1')$ obtained contracting
the edge $e_1$ in $\Gamma_0$ and $\Gamma_1$ respectively.
By construction the contractions $\Gamma_0'\to\Gamma_1'$ and
$\Gamma_0'\to\Gamma_2$ still respect the hypothesis, and $\#V(\Gamma_1')=\#V(\Gamma_1)-1$.
Therefore by induction the automorphism $\msf a'$,
which is an element of $G(\Gamma_0',M')$, is a sum
$$\msf a'=\msf a_1' +\msf a_2,$$
with $\msf a_1'\in G(\Gamma_1')\subset G(\Gamma_1)$ and $\msf a_2\in G(\Gamma_2)$.
Finally $\msf a=k\cdot\cut_{\Gamma_1}(e;T_1)+\msf a_1'+\msf a_2$, then it is in $G(\Gamma_1)+G(\Gamma_2)$.
This is a contradiction 
and so the lemma is proved.\fine

In what follows we will find, for some small prime values of $\ell$,
a description of $J_{g,\ell}^k$ by our stratification. Before
starting we point out that our analysis will focus in the cases
$J_{g,\ell}^0$ and~$J_{g,\ell}^1$.
\begin{prop}\label{prop_k}
If $\ell$ is prime, we have
a natural stack isomorphism $\rr_{g,\ell}^1\cong\rr_{g,\ell}^k$ for every
$k$ between $1$ and $\ell-1$.
\end{prop}
\proof Consider a scheme $S$ and a triple $(\ssC\to S,\ssL,\phi)$ in $\rr_{g,\ell}^1(S)$.
Consider the map sending it to
$$(\ssC\to S,\ssL^{\xx k},\phi^{\xx k})\in \rr_{g,\ell}^k(S).$$
As $\ell$ is prime and $k\not\equiv 0\mod\ell$, this morphism
has a canonical inverse, so we obtained an isomorphism of categories.\fine

\subsection{The locus $J_{g,\ell}^k$ for $\ell=2$}
In this case the $J$-locus is always empty. Indeed, every 
automorphisms in $\underline\Aut_C(\ssC,\ssL,\phi)\slash \QR$
 must have a support of cardinality at least $2$,
 but for every edge in the support, a ghost $\msf a$ has
a contribution of $1\slash 2$ to its age. Hence there are 
no junior automorphisms in $\underline\Aut_C(\ssC,\ssL,\phi)\slash\QR$.
This result was already obtained by Farkas and Ludwig
for the Prym space $\overline\R_{g,2}^0$  in \cite{farlud10}, and by Ludwig for $\overline\R_{g,2}^1$
in \cite{lud10}.
\newline

\subsection{The locus $J_{g,\ell}^k$ for $\ell=3$}
The process of finding the $J_{g,\ell}^k$ decomposition, for
specific values of $\ell$ and $k$, will always follows three steps.
\newline

\cor{Step 1.} We identify at first the graphs which can
support a junior automorphism in $\underline\Aut(\ssC,\ssL,\phi)\slash\QR$,
\cor{i.e.}~those graphs with $\#E<\ell$ and no separating edges.
If $\Gamma$ is one of these graphs, we identify the
$\Z\slash\ell$-valued automorphisms supported on $\Gamma$,
\cor{i.e.}~the elements of $\bigoplus_{e\in E(\Gamma)}\Z\slash\ell$
which are non-trivial on every edge and are junior. These
automorphisms are the junior elements in
 $\Aut_C(\ssC)$ for an $\ell$-twisted curve
whose dual graph is $\Gamma$.

If $\ell=3$ there is only one junior automorphism
which can be supported on a $\Z\slash 3$-valued
decorated graph, the one represented
in the image below and supported on $\Gamma_{(2,2)}$.
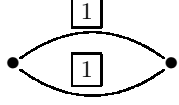
\begin{figure}[h]
 \begin{picture}(400,25)(-40,20)
 
   \put(130,30){\xymatrix{*{\bullet} \ar@{-}@/^1pc/[rr] \ar@{-}@/_1pc/[rr]  
   & &*{\bullet}}}
   \put(155,47){\framebox{\footnotesize{1}}}
    \put(155,25){\framebox{\footnotesize{1}}}
 
  \end{picture}
\end{figure}

\cor{Step 2.} For each one of these junior automorphisms,
we search
for multiplicity cochains that respect the lift condition of
Theorem \ref{lift_condition} on the automorphisms above.
In this case the only possibilities are the following
cochains

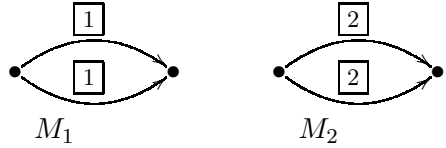
\begin{figure}[h]
 \begin{picture}(400,35)(-40,13)
 
   \put(60,30){\xymatrix{*{\bullet} \ar@/^1pc/[rr] \ar@/_1pc/[rr]  
   & &*{\bullet}}}
   \put(85,47){\framebox{\footnotesize{1}}}
    \put(85,25){\framebox{\footnotesize{1}}}
    \put(70,5){$M_1$}
    
       \put(160,30){\xymatrix{*{\bullet} \ar@/^1pc/[rr] \ar@/_1pc/[rr]  
   & &*{\bullet}}}
   \put(185,47){\framebox{\footnotesize{2}}}
    \put(185,25){\framebox{\footnotesize{2}}}
    \put(170,5){$M_2$}
 
  \end{picture}
\end{figure}

In fact, the two decorated graphs $(\Gamma_{(2,2)},M_1)$ and $(\Gamma_{(2,2)},M_2)$
are isomorphic by the isomorphism inverting the two vertices.
Thus for $\ell=3$, there is only one class of decorated graphs admitting
junior automorphisms.

\cor{Step 3.}
By Proposition \ref{prop_mdeg}, there is an additional condition that
the decorated graph must satisfy: the cochain $\partial M_1$ have
to be the multidegree cochain of $\omega_\ssC^{\xx k}$. Equivalently
\begin{equation}\label{multideg}
\sum_{e_+=v} M_1(e)\equiv \deg\left.\omega\right|_v^{\xx k}\equiv k\cdot(2g_v-2+N_v)\ \mod\ell \ \ \ \ \forall v\in V(\Gamma),
\end{equation}
where $g_v$ is the genus of the component correspondent to vertex $v$, and $N_v$ is the
degree of this vertex, \cor{i.e.}~the number of edges touching it.

By Proposition \ref{prop_k} we can focus in cases
$k=0$ and $1$. If $k=1$ the condition of (\ref{multideg}) is empty,
because $2$ and $3$ are coprime and there always exists a sequence of $g_v$
satisfying the equality. Then we have
$$J_{g,3}^1=\overline\s_{(\Gamma_{(2,2)},M_1)}.$$

In case $k=0$, by (\ref{multideg}) we have $\sum_{e_+=v}M_1(e)\equiv 0\mod 3$
for both vertices, but this condition is not satisfied by $\left(\Gamma_{(2,2)},M_1\right)$,
then 
$$J_{g,3}^0=\varnothing.$$

\subsection{The locus $J_{g,\ell}^k$ for $\ell=5$}	

\cor{Step 1 and 2.}
For $\ell=5$, 
every graph $(\Gamma,M)$ such that there exists a junior automorphism $\msf a\in G(\Gamma,M)$,
contracts on a vine stratum. This is a consequence of the following lemma.

\begin{lemma}\label{lem_config}
If $\ell$ is a prime number,
consider a decorated graph $(\Gamma,M)$ such that there exists
a vertex $v_1\in V(\Gamma)$ connected with exactly two vertices $v_2,v_3$
and such that between $v_1$ and $v_2$ there is only one edge called~$e$
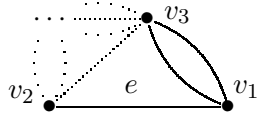
\begin{figure}[h]
 \begin{picture}(400,25)(0,0)
 
   \put(180,30){\xymatrix{{\cdots}\ar@{.}@/^/[r] \ar@{.}[r]  \ar@{.}@/_/[r] \ar@{.}@/^/[d] \ar@{.}@/_/[d] &*{\bullet} \ar@{-}@/^/[rd] \ar@{-}@/_/[rd] &\\
  *{\bullet}\ar@{-}[rr] \ar@{.}[ru]& &*{\bullet}}}
   \put(174,0){$v_2$}
    \put(259,2){$v_1$}
    \put(218,2){$e$}
    \put(233,30){$v_3$}
    
  \end{picture}
\end{figure}

If $\s_{(\Gamma,M)}$ is a junior stratum, then there exists a non-trivial graph contraction $(\Gamma,M)\to(\Gamma_1,M_1)$
such that $\s_{(\Gamma_1,M_1)}$ is also junior.
\end{lemma}
\begin{rmk}
This lemma permits one to simplify the analysis of junior strata.
Every stratum labeled with a graph containing the configuration above, can
be ignored in the analysis. Indeed, if $\s_{(\Gamma,M)}$ is a subset of the $J$-locus,
there exists a decorated graph $(\Gamma_1,M_1)$ with less vertices such that
$\s_{(\Gamma,M)}\subset\overline\s_{(\Gamma_1,M_1)}\subset J_{g,\ell}^k$.
\end{rmk}
\proof[Proof of the Lemma]
Consider $\msf a\in G(\Gamma,M)$ such that $\age \msf a<1$. If $\msf a$ is not supported
on $\s_{(\Gamma,M)}$, we contract one edge where $\msf a$ acts trivially and the lemma is proved.
Thus we suppose $\msf a$ supported on $\s_{(\Gamma,M)}$.
We call $e'$ one of the edges connecting $v_1$ and $v_3$. We consider a spanning tree $T$ of $\Gamma$
passing though $e'$ and not passing through $e$. Then we call $\Gamma_1$ and $\Gamma_2$
the two contractions of $\Gamma$ obtained contracting respectively $e'$ and $E_T\backslash\{e'\}$.
By Lemma \ref{lem_cover}, we have
$$G(\Gamma,M)=G(\Gamma_1,M_1)\oplus G(\Gamma_2,M_2).$$
Therefore we use Theorem \ref{teo_levage} to obtain
$$\age\msf a-\#E(\Gamma)\geq\left(\age\s_{(\Gamma_1,M_1)}-\#E(\Gamma_1)\right)+
\left(\age\s_{(\Gamma_2,M_2)}-\#E(\Gamma_2)\right).$$
As $\#E(\Gamma)=\#E(\Gamma_1)+\#E(\Gamma_2)-1$ by construction,
$\msf a$ junior implies $\s_{(\Gamma_1,M_1)}$ or $\s_{(\Gamma_2,M_2)}$ to
be junior.\fine

The configuration
of Lemma \ref{lem_config} appears
in every non-vine graph with less than $5$ edges.
As a consequence
the reduction to vine strata follows.

To summarize these vine strata we introduce a new notation.
Consider a $k$-vine graph, with vertices $v_1$ and $v_2$
and edges $e_1,\dots,e_k$ all taken oriented from $v_1$ to $v_2$.
If the multiplicity index $M$ on the graph take values
$M(e_i)=m_i$, $m_i\in\Z\slash\ell$, we 
note this decorated graph $(m_1,m_2,\dots,m_k)$.
For example the graphs appeared in the precedent paragraph
are noted $(1,1)$ and $(2,2)$ (and are isomorphic).
We can now state the ten classes of vine graphs which support a junior
ghost for $\ell=5$,
$$(1,1),\ (2,2),\ (1,2),\ (1,3),\ (1,1,1),\ (2,2,2),\ (1,1,3),\ (2,2,1),\ (1,1,1,1),\ (2,2,2,2).$$

\cor{Step 3.}
For $k=1$, equation (\ref{multideg}) is always
respected for some genus labellings of the graph. Therefore we have
$$J_{g,5}^1=\overline\s_{(1,1)}\cup\overline\s_{(2,2)}\cup\overline\s_{(1,2)}\cup
\overline\s_{(1,3)}\cup \overline\s_{(1,1,1)}\cup\overline\s_{(2,2,2)}\cup\overline\s_{(1,1,3)}\cup\overline\s_{(2,2,1)}\cup \overline\s_{(1,1,1,1)}\cup\overline\s_{(2,2,2,2)}.$$
If $k=0$ the equation is satisfied by two vine graphs, and we obtain the following
$$J_{g,5}^0=\overline\s_{(1,1,3)}\cup\overline\s_{(2,2,1)}.$$
In particular, this result fills the hole in Chiodo and Farkas analysis in \cite{chiofar12}.
They proved that the $J$-locus in the space $\overline\R_{g,\ell}^0$
is empty for $\ell\leq 6$ and $\ell\neq 5$.

\subsection{The locus $J_{g,\ell}^k$ for $\ell=7$}
\cor{Step 1 and 2.} Using Lemma \ref{lem_config}, we observe
that for $\ell=7$ there are two kinds of graphs admitting junior automorphisms. 
The first kind is the usual vine graph, but we can also have a $3$-cycle graph
such that every pair of vertices is connected by two edges. In this second case, the
only possible automorphism with age lower than $1$ takes value $1$ on 
every edge. As a consequence the possible decorations are like in figure below.

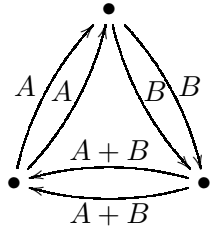
\begin{figure}[h]
 \begin{picture}(400,72)(0,10)
 
   \put(180,80){\xymatrix{ & {\bullet}\ar@/^/[rdd] \ar@/_/[rdd] \\
   \\
   {\bullet}\ar@/^/[ruu] \ar@/_/[ruu] & & {\bullet}\ar@/^/[ll] \ar@/_/[ll]}
   }
   \put(186,47){$A$}
   \put(200,45){$A$}
    \put(235,45){$B$}
    \put(248,47){$B$}
    \put(207,23){$A+B$}
    \put(207,-1){$A+B$}
    
  \end{picture}
\end{figure}
We define two sets of decorated graphs
$$V_7:=\{\m{vine decorated graphs admitting junior automorphism}\}\slash\cong$$
$$C_7:=\{\m{graphs decorated as in the figure}\}\slash\cong.$$

\cor{Step 3.}
If $k=1$, condition (\ref{multideg}) is always verified by some genus labeling.
Therefore, we have
$$J_{g,7}^1=\bigcup_{(\Gamma,M)\in V_7}\overline\s_{(\Gamma,M)}\cup\bigcup_{(\Gamma,M)\in C_7}\overline\s_{(\Gamma,M)}.$$

If $k=0$, we call $V_7'$ the subset of $V_7$ of decorated graphs
respecting equation (\ref{multideg}). Every graph in $C_7$ does not respect
the equation. Indeed, we must have $4A+2B\equiv0\mod7$ and $2B-2A\equiv 0\mod7$,
therefore $A\equiv B\equiv 0$ which is not allowed. Finally, we have
$$J_{g,7}^0=\bigcup_{(\Gamma,M)\in V_7'}\overline\s_{(\Gamma,M)}.$$
In other words, the $J$-locus of $\overline\R_{g,7}^0$ is a union
of vine strata.

\bibliography{mybiblio}{}

\begin{thebibliography}{10}

\bibitem{acv03}
Dan Abramovich, Alessio Corti, and Angelo Vistoli.
\newblock Twisted bundles and admissible covers.
\newblock {\em Comm. Algebra}, 31(8):3547--3618, 2003.
\newblock Special issue in honor of Steven L. Kleiman.

\bibitem{abrajar03}
Dan Abramovich and Tyler~J. Jarvis.
\newblock Moduli of twisted spin curves.
\newblock {\em Proc. Amer. Math. Soc.}, 131(3):685--699, 2003.

\bibitem{abravis02}
Dan Abramovich and Angelo Vistoli.
\newblock Compactifying the space of stable maps.
\newblock {\em J. Amer. Math. Soc.}, 15(1):27--75, 2002.

\bibitem{arbacor87}
Enrico Arbarello and Maurizio Cornalba.
\newblock The {P}icard groups of the moduli spaces of curves.
\newblock {\em Topology}, 26(2):153--171, 1987.

\bibitem{biggs93}
Norman Biggs.
\newblock {\em Algebraic graph theory}.
\newblock Cambridge Mathematical Library. Cambridge University Press,
  Cambridge, second edition, 1993.

\bibitem{chio08}
Alessandro Chiodo.
\newblock Stable twisted curves and their {$r$}-spin structures.
\newblock {\em Ann. Inst. Fourier (Grenoble)}, 58(5):1635--1689, 2008.

\bibitem{cefs13}
Alessandro Chiodo, David Eisenbud, Gavril Farkas, and Frank-Olaf Schreyer.
\newblock Syzygies of torsion bundles and the geometry of the level {$\ell$}
  modular variety over {$\overline{\mc{M}}_g$}.
\newblock {\em Invent. Math.}, 194(1):73--118, 2013.

\bibitem{chiofar12}
Alessandro Chiodo and Gavril Farkas.
\newblock Singularities of the moduli space of level curves.
\newblock {\em to appear in J.~Eur. Math. Soc.}, 2015.
\newblock arXiv:1205.0201.

\bibitem{far00}
Gavril Farkas.
\newblock The geometry of the moduli space of curves of genus 23.
\newblock {\em Math. Ann.}, 318(1):43--65, 2000.

\bibitem{farlud10}
Gavril Farkas and Katharina Ludwig.
\newblock The {K}odaira dimension of the moduli space of {P}rym varieties.
\newblock {\em J. Eur. Math. Soc. (JEMS)}, 12(3):755--795, 2010.

\bibitem{farver10}
Gavril Farkas and Alessandro Verra.
\newblock The geometry of the moduli space of odd spin curves.
\newblock {\em Ann. of Math. (2)}, 180(3):927--970, 2014.

\bibitem{harmum82}
Joe Harris and David Mumford.
\newblock On the {K}odaira dimension of the moduli space of curves.
\newblock {\em Invent. Math.}, 67(1):23--88, 1982.
\newblock With an appendix by William Fulton.

\bibitem{lud10}
Katharina Ludwig.
\newblock On the geometry of the moduli space of spin curves.
\newblock {\em J. Algebraic Geom.}, 19(1):133--171, 2010.

\bibitem{mestrarama85}
N.~Mestrano and S.~Ramanan.
\newblock Poincar\'e bundles for families of curves.
\newblock {\em J. Reine Angew. Math.}, 362:169--178, 1985.

\bibitem{prill67}
David Prill.
\newblock Local classification of quotients of complex manifolds by
  discontinuous groups.
\newblock {\em Duke Math. J.}, 34:375--386, 1967.

\bibitem{reid80}
Miles Reid.
\newblock Canonical {$3$}-folds.
\newblock In {\em Journ\'ees de {G}\'eometrie {A}lg\'ebrique d'{A}ngers,
  {J}uillet 1979/{A}lgebraic {G}eometry, {A}ngers, 1979}, pages 273--310.
  Sijthoff \& Noordhoff, Alphen aan den Rijn, 1980.

\end{thebibliography}
\bibliographystyle{plain}

\end{document}